                    \def\version{March 6, 2002}                          %
\documentclass[reqno,11pt]{amsart}

\usepackage{amsmath,amsfonts}
\usepackage{amssymb}
\usepackage[dvips]{graphicx}

\newfam\Bbbfam
\font\tenBbb=msbm10 \font\sevenBbb=msbm7 \font\fiveBbb=msbm5
\textfont\Bbbfam=\tenBbb \scriptfont\Bbbfam=\sevenBbb
\scriptscriptfont\Bbbfam=\fiveBbb

\newcommand{\R}     {{\mathbb{R}}}

\newcommand{\N}     {{\mathbb{N}}}
\newcommand{\Q}     {{\mathbb{Q}}}

\renewcommand{\P}   {{\mathbb{P}}}

\newcommand{\E}     {{\mathbb{E}}}

\newcommand{\Ai}    {{\rm{ Ai}}}
\newcommand{\Bi}    {{\rm{ Bi}}}

\def\1{{\mathchoice {1\mskip-4mu\mathrm l}
{1\mskip-4mu\mathrm l} {1\mskip-4.5mu\mathrm l}
{1\mskip-5mu\mathrm l}}}

\def\comment#1{}
\newtheoremstyle{thm}{2ex}{2ex}{\itshape\rmfamily}{}
{\bfseries\rmfamily}{}{1.7ex}{}

\newtheoremstyle{rem}{1.3ex}{1.3ex}{\rmfamily}{}
{\itshape\rmfamily}{}{1.5ex}{}

\newenvironment{proofsect}[1]
{\vskip0.1cm\noindent{\bf #1}}{\vspace{0.15cm}}

\newtheorem{theorem}{Theorem}[section]
\newtheorem{lemma}[theorem]{Lemma}
\newtheorem{prop}[theorem] {Proposition}
\newtheorem{cor}[theorem]  {Corollary}

\newtheorem{step}{STEP}

\newcommand{\en}       {\end{equation}}
\newcommand{\eq}       {\begin{equation}}

\newcommand{\eqry}     {\begin{eqnarray}}
\newcommand{\enqry}    {\end{eqnarray}}
\newcommand{\eqarray}  {\begin{eqnarray}}
\newcommand{\enarray}  {\end{eqnarray}}
\newcommand{\eqarraystar} {\begin{eqnarray*}}
\newcommand{\enarraystar} {\end{eqnarray*}}
\newcommand{\bel}      {\begin{lemma}}
\newcommand{\el}       {\end{lemma}}

\newcommand{\bec}      {\begin{cor}}
\newcommand{\ec}       {\end{cor}}
\newcommand{\bes}      {\begin{step}}
\newcommand{\es}       {\end{step}}
\newcommand{\bea}      {\begin{array}}
\newcommand{\ea}       {\end{array}}
\newcommand{\bpr}      {\begin{proof}}
\newcommand{\epr}      {\end{proof}}

\newcommand{\da}       {\downarrow}

\renewcommand{\section}{\secdef\sct\sect}
\newcommand{\sct}[2][default]{\refstepcounter{section}
\vspace{0.5cm} \setcounter{equation}{0}
\centerline{\scshape \arabic{section}.\ #1} \vspace{0.3cm}}
\newcommand{\sect}[1]{\vspace{0.5cm}
\centerline{\large\scshape #1} \vspace{0.3cm}}

\renewcommand{\subsection}{\secdef\subsct\sbsect}
\newcommand{\subsct}[2][default]{\refstepcounter{subsection}
\nopagebreak \vspace{0.5\baselineskip} {\flushleft\bf
\arabic{section}.\arabic{subsection}~\bf #1 } \nopagebreak}
\newcommand{\sbsect}[1]{\vspace{0.1cm}\noindent
{\bf #1}\vspace{0.1cm}}

\renewcommand{\subsubsection}{\secdef\subsubsect\sbsbsect}
\newcommand{\subsubsect}[2][default]{\refstepcounter{subsubsection}
\nopagebreak \vspace{0.1\baselineskip} {\flushleft\sffamily\slshape
\arabic{section}.\arabic{subsection}.\arabic{subsubsection}
\sffamily #1\/.} }
\newcommand{\sbsbsect}[1]{\vspace{0.1cm}\noindent
{\bf #1}\ }

\newcommand{\ef}             {{\rm e}}
\newcommand{\eps}            {\varepsilon}

\renewcommand{\d}            {{\rm d}}

\newcommand{\Ccal}   {{\mathcal C }}

\newcommand{\Ecal}   {{\mathcal E }}

\newcommand{\Kcal}   {{\mathcal K }}

\newcommand{\Ocal}   {{\mathcal O }}

\newcommand{\eve}    {{{\Ecal}}}
\newcommand{\evf}    {{{\Ecal}}}

\newcommand{\smallsup}[1] {{\scriptscriptstyle{({#1}})}}

\begin{document}

\title[Large deviations for the one-dimensional Edwards model]
{\large Large deviations for the\\
\vskip 0.3cm
one-dimensional Edwards model}

\author[Remco van der Hofstad, Frank den Hollander, Wolfgang K\"onig]{}
\maketitle
\thispagestyle{empty}
\vspace{1cm}

\centerline{\small\version}

\vspace{0.5cm}

\begin{center}
Remco van der Hofstad
\footnote{Department of Applied Mathematics,
Delft University of Technology, Mekelweg 4,
2628 CD Delft, The Netherlands.}
\footnote
{Present address: Department of Mathematics and Computer Science,
Eindhoven University of Technology, P.O.\ Box  513,
5600 MB Eindhoven, The Netherlands.
{\tt rhofstad@win.tue.nl}.}\\
Frank den Hollander
\footnote{EURANDOM, P.O.\ Box 513,
5600 MB Eindhoven, The Netherlands.
{\tt denhollander@eurandom.tue.nl}}\\
Wolfgang K\"onig \footnote{Institut f\"ur Mathematik, TU Berlin,
Stra{\ss}e des 17.\ Juni 136, D-10623 Berlin, Germany. {\tt
koenig@math.tu-berlin.de}}

\end{center}

\vspace{2cm}

{\small
\begin{quote}
{\bf Abstract:} In this paper we prove a large deviation principle
for the empirical drift of a one-dimensional Brownian motion with
self-repellence called the Edwards model. Our results extend
earlier work in which a law of large numbers, respectively, a
central limit theorem were derived. In the Edwards model a path of
length $T$ receives a penalty $e^{-\beta H_T}$, where $ H_T$ is
the self-intersection local time of the path and
$\beta\in(0,\infty)$ is a parameter called the strength of
self-repellence. We identify the rate function in the large
deviation principle for the endpoint of the path as $\beta^{\frac
23} I(\beta^{-\frac 13}\cdot)$, with $I(\cdot)$ given in terms of
the principal eigenvalues of a one-parameter family of
Sturm-Liouville operators. We show that there exist numbers
$0<b^{**}<b^*<\infty$ such that: (1) $I$ is linearly decreasing on
$[0,b^{**}]$; (2) $I$ is real-analytic and strictly convex on
$(b^{**},\infty)$; (3) $I$ is continuously differentiable at
$b^{**}$; (4) $I$ has a unique zero at $b^*$. (The latter fact
identifies $b^*$ as the asymptotic drift of the endpoint.) The
critical drift $b^{**}$ is associated with a crossover in the
optimal strategy of the path: for $b\geq b^{**}$ the path assumes
local drift $b$ during the full time $T$, while for $0\leq
b<b^{**}$ it assumes local drift $b^{**}$ during time
$\frac{b^{**}+b}{2b^{**}}T$ and local drift $-b^{**}$ during the
remaining time $\frac{b^{**}-b}{2b^{**}}T$. Thus, in the second
regime the path makes an overshoot of size $\frac{b^{**}-b}{2}T$
in order to reduce its intersection local time.
\end{quote}
\normalsize}

\vfill


 \eject

\setcounter{section}{0}

%%%%%%%%%%%%%%%% SECTION 1 %%%%%%%%%%%%%%%%%%%%%%%%%%%%%%%%%%%%%%%%%%%

\section{Introduction and main results}
\label{EDW}

\subsection{The Edwards model}
\label{sec-EM}

Let $B=(B_t)_{t \geq 0}$ be standard Brownian motion on $\R$
starting at the origin ($B_0=0$). Let $ P$ be the Wiener measure
and let $ E$ be expectation with respect to $ P$. For $T>0$ and
$\beta\in(0,\infty)$, define a probability law $\Q_T^{\beta}$ on
paths of length $T$ by setting
    \begin{equation}
    \frac{\d\Q_T^\beta}{\d P}[\,\cdot\,]
    =\frac{1}{ Z_T^\beta}e^{-\beta H_T[\,\cdot\,]},
    \qquad  Z_T^\beta= E(e^{-\beta H_T}),
    \label{Ed mod}
    \end{equation}
where
    \begin{equation}
    \label{HT}
     H_T\left[(B_t)_{t\in[0,T]}\right]
    =\int_0^T{\rm d}u \int_0^T {\rm d}v\,\,\delta(B_u-B_v)
    =\int_{\R} {\rm d}x\,\,L(T,x)^2
    \end{equation}
is the Brownian intersection local time up to time $T$. The first
expression in (\ref{HT}) is formal only. In the second expression
the Brownian local times $L(T,x)$, $x\in\R$, appear. The law $\Q_T^\beta$
is called the $T$-polymer measure with strength of self-repellence $\beta$.
The Brownian scaling property implies that
    \begin{equation}
    \label{scal}
    \Q_T^\beta\left((B_t)_{t\in[0,T]}\in \cdot\,\right)
    = \Q_{\beta^{\frac23}T}^1
    \left((\beta^{-\frac13}B_{\beta^{\frac23}t})_{t\in[0,T]}\in\cdot\,\right).
    \end{equation}

It is known that under the law $\Q_T^\beta$ the endpoint $B_T$ satisfies the
following central limit theorem:

    \medskip
    \begin{theorem}[Central limit theorem]
    \label{thm-LDEM*}
    There are numbers $a^*,b^*,c^*\in(0,\infty)$ such that for any $\beta\in(0,\infty)$:
    \begin{enumerate}
    \item[(i)]
    Under the law $\Q_T^\beta$, the distribution of the scaled endpoint
    $(|B_T|-b^*\beta^{\frac 13}T)/c^*\sqrt T$ converges weakly to the
    standard normal distribution.
    \item[(ii)]
    $\lim_{T\to\infty}\frac 1T\log  Z_T^\beta=-a^*\beta^{\frac 23}$.
    \end{enumerate}
    \end{theorem}

\medskip\noindent
Theorem~\ref{thm-LDEM*} is contained in \cite[Theorem 2 and Proposition 1]{HHK1}.
For the identification of $a^*, b^*, c^*$, see \eqref{constdefs} below. Bounds on
these numbers appeared in \cite[Theorem~3]{Hconst}. The numerical values are:
$a^*\approx 2.19$, $b^*\approx 1.11$, $c^*\approx 0.63$. The law of large numbers
corresponding to Theorem~\ref{thm-LDEM*}(i) was first obtained by Westwater \cite{We3} (see also
\cite[Section~0.6]{HH}).

The main object of interest in the present paper is the rate function $
I_\beta$ defined by
\footnote{In fact, $ I_\beta$ differs by a constant from what is usually
called a rate function. This constant is $\lim_{T\to\infty}\frac{1}{T}\log
 Z_T^\beta=-a^*\beta^{\frac 23}$ (see Theorem \ref{thm-LDEM*}(ii)).
Hence, $I_\beta -a^*\beta^{\frac 23}$ is the true rate function.}
    \begin{equation}
    \label{rateED}
    \begin{aligned}
    - I_\beta(b) &= \lim_{T\to\infty}\frac 1T\log
    E\Bigl(e^{-\beta H_T} \1_{\{B_T \approx bT\}}\Bigr)\\
    &= \lim\limits_{T\to\infty}\frac 1T\log
    \left\{ Z^\beta_T \Q_T^\beta(B_T\approx bT)\right\},
    \qquad b\in\R,
    \end{aligned}
    \end{equation}
where $B_T \approx bT$ is an abbreviation for $|B_T-bT|\leq \gamma_T$ for some
$\gamma_T>0$ such that $\gamma_T/T\to 0$ and $\gamma_T/\sqrt{T}\to\infty$ as
$T\to\infty$. (We will see that the limit in \eqref{rateED} does not depend
on the choice of $\gamma_T$.) It is clear from \eqref{scal} that
    \begin{equation}
    \label{Iscal}
    \beta^{-\frac 23} I_\beta(\beta^{\frac 13}b)
    = I_1(b), \qquad b\geq 0,
    \end{equation}
provided the limit in \eqref{rateED} exists for $\beta=1$ and $b\geq 0$. Moreover,
    \begin{equation}
     I_\beta(b)= I_\beta(-b), \qquad b\leq 0.
    \end{equation}
Therefore, we may restrict ourselves to $\beta=1$ and $b\geq 0$. In the following
we write $ I =  I_1$.

\subsection{Main results}
\label{sec-EMtheorems}

Our first main result says that $ I$ exists and has the shape exhibited
in Fig.~1. \footnote{In \cite[Corollary~2.6 and Remark~2.7]{MS}
it was proved that $\lim_{T\to\infty}\frac 1T\log  E(e^{-H_T}\,
|\,B_T=0)=-a^{**}$, which essentially gives the existence of $I(0)$ with
value $a^{**}$. Furthermore, the existence of $I(b^*)$ with
value $a^{*}$ follows from our earlier work \cite[Proposition~1]{HHK1}.}

    \medskip
    \begin{theorem}[Large deviations]
    \label{thm-LDEM}
    Let $\beta=1$.
    \begin{enumerate}
    \item[(i)]
    For any $b\geq 0$, the limit $I(b)$ in \eqref{rateED} exists and
    is finite.
    \item[(ii)]
    $ I$ is continuous and convex on $[0,\infty)$, and continuously
    differentiable on $(0,\infty)$.
    \item[(iii)]
    There are numbers $a^{**}\in(a^*,\infty)$, $b^{**}\in(0,b^*)$ and $\rho(a^{**})
    \in(0,\infty)$ such that $ I(0)=a^{**}$, $ I$ is linearly decreasing
    on $[0,b^{**}]$ with slope $-\rho(a^{**})$, is real-analytic and strictly convex
    on $(b^{**},\infty)$, and attains its unique minimum at $b^*$ with $
    I(b^*)=a^*$ and $ I\,''(b^*)=1/c^{*2}$.
    \item[(iv)] $I(b)=\frac12 b^2+\Ocal(b^{-1})$ as $b\to\infty$.
    \end{enumerate}
    \end{theorem}

%%%%%%%%%%%%%%%%%%% BEGIN FIGURE %%%%%%%%%%%%%%%%%%%%%%%%%%%%%%%%%%%%%%%%%%%%%%%%%

\begin{center}
\vskip 1truecm

\setlength{\unitlength}{0.7cm}

\begin{picture}(14,8)(-1,-3)%(12,8)(-3,-3)

  \put(0,-2){\line(12,0){12}}
  \put(0,-2){\line(0,7){7}}
  \put(-.1,-3){$0$}
  {\thicklines
   \qbezier(0,3)(2,2.5)(4,2)
   \qbezier(4,2)(6,1.55)(7,1.85)
   \qbezier(7,1.85)(8,2.2)(9,5)
  }
  \qbezier[40](4,-2)(4,0)(4,2)
  \qbezier[40](6,-2)(6,-.25)(6,1.75)
  \qbezier[70](0,1.7)(4,1.7)(8,1.7)
  \qbezier[70](4,2)(6,1.5)(7,1.25)

  \put(4,2){\circle*{.15}}
  \put(6,1.74){\circle*{.15}}
  \put(3.8,-3){$b^{**}$}
  \put(5.9,-3){$b^*$}
  \put(12.5,-2.15){$b$}
  \put(-1,1.62){$a^*$}
  \put(-1,2.9){$a^{**}$}
  \put(-.25,5.5){$I(b)$}
  \put(7.4,1){slope $-\rho(a^{**})$}

  \end{picture}
  \end{center}

\begin{center}
\small Fig.\ 1. Qualitative picture of $b\mapsto I(b)$.
\end{center}

\vskip 0.5truecm

%%%%%%%%%%%%%%%%%%% END FIGURE %%%%%%%%%%%%%%%%%%%%%%%%%%%%%%%%%%%%%%%%%%%%%%%%
\medskip

The linear piece of the rate function has the following intuitive interpretation.
If $b\geq b^{**}$, then the best strategy for the path to realize the large deviation
event $\{B_T\approx bT\}$ is to assume local drift $b$ during time $T$. In particular,
the path makes no overshoot on scale $T$, and this leads to the real-analyticity and
strict convexity of $I$ on $(b^{**},\infty)$. On the other hand, if $0\leq b<b^{**}$,
then this strategy is too expensive, since too small a drift leads to too large an
intersection local time. Therefore the best strategy now is to assume local drift
$b^{**}$ during time $\frac{b^{**}+b}{2b^{**}}T$ and local drift $-b^{**}$ during
the remaining time $\frac{b^{**}-b}{2b^{**}}T$. In particular, the path makes
an overshoot on scale $T$, namely, $\frac{b^{**}-b}{2}T$, and this leads to the
linearity of $I$ on $[0,b^{**}]$. At the critical drift $b=b^{**}$, $I$ is
continuously differentiable.

For $b\to\infty$, $I(b)$ is determined by the Gaussian tail of $B_T$ because the
intersection local time $H_T$ vanishes.

For the identification of $a^{**},b^{**},\rho(a^{**})$, see \eqref{b**def} below.
The numerical values are: $a^{**}\approx~2.95$, $b^{**}\approx~0.85$, $\rho(a^{**})
\approx~0.78$. These estimates can be obtained with the help of the method
in \cite{Hconst}.

There is an intimate connection between the rate function $ I$ and the two
moment generating functions $\Lambda^{+},\Lambda^- \colon\R\to\R$ given by
    \begin{equation}
    \Lambda^+(\mu) = \lim_{T\to\infty}\frac 1T\log
     E\bigl(e^{- H_T}e^{\mu B_T}\1_{\{B_T\geq 0\}}\bigr),\label{Lambda+def}
    \end{equation}
and the same formula for $\Lambda^-(\mu)$ with $\1_{\{B_T\geq 0\}}$ replaced
by $\1_{\{B_T\leq 0\}}$. Obviously, $\Lambda^+(-\mu) = \Lambda^-(\mu)$ for any $\mu\in\R$,
provided one limit exists.

Our second main result says that $\Lambda^+$ exists and has the shape exhibited in Fig.~2,
and that its Legendre transform is equal to $ I$ on $[0,\infty)$.

    \medskip
    \begin{theorem}[Exponential moments]
    \label{expmom}
    Let $\beta=1$.
    \begin{enumerate}
    \item[(i)]
    For any $\mu\in\R$, the limit $\Lambda^+(\mu)$ in \eqref{Lambda+def} exists and
    is finite.
    \item[(ii)]
    $\Lambda^+$ equals $-a^{**}$ on $(-\infty,-\rho(a^{**})]$, is real-analytic and
    strictly convex on $(-\rho(a^{**}),\infty)$, and satisfies $\lim_{\mu\da -\rho(a^{**})}
    (\Lambda^{+})'(\mu)=b^{**}$.
    \item[(iii)]
    $\Lambda^{+}(\mu)=\frac12 \mu^2 +\Ocal(\mu^{-1})$ as $\mu\to\infty$.
    \item[(iv)]
    The restriction of  $ I$ to $[0,\infty)$ is the Legendre transform of $\Lambda^+$,
    i.e.,
    \begin{equation}
    \label{LegTra}
    I(b)=\max_{\mu\in \R}\bigl[b\mu-\Lambda^+(\mu)\bigr], \qquad b \geq 0.
    \end{equation}
    \end{enumerate}
    \end{theorem}

%%%%%%%%%%%%%%%%%%% BEGIN FIGURE %%%%%%%%%%%%%%%%%%%%%%%%%%%%%%%%%%%%%%%%%%%%%%%%%

%\vskip 1truecm

\setlength{\unitlength}{0.7cm}

\begin{picture}(8,12)(-11,-6)

  \put(0,1){\line(5,0){5}}
  \put(0,1){\line(-5,0){5}}
  \put(0,0){\line(0,4){4}}
  \put(0,0){\line(0,-4){4}}
  \put(-.4,.6){$0$}
  {\thicklines
   \qbezier(-5,-3)(-4,-3)(-2,-3)
   \qbezier(-2,-3)(0,-2.8)(4,3.5)
  }
  \qbezier[40](-2,1)(-2,-.5)(-2,-3)
  \qbezier[15](-2,-3)(-1,-3)(-.1,-3)
  \qbezier[70](-4,-3.4)(-2,-3)(2,-2.2)
  \qbezier[70](-4,-5.4)(0,-1.8)(2,0)

  \put(-2,-3){\circle*{.15}}
  \put(0,-1.8){\circle*{.15}}
  \put(-2.4,1.3){$-\rho(a^{**})$}
  \put(.3,-1.9){$-a^*$}
  \put(5.5,.9){$\mu$}
  \put(.3,-3.1){$-a^{**}$}
  \put(-.35,4.5){$\Lambda^+(\mu)$}
  \put(2.4,-2.05){slope $b^{**}$}
  \put(2.4,-.2){slope $b^*$}

  \end{picture}
\vskip -1truecm
\begin{center}
\small
Fig.\ 2. Qualitative picture of $\mu\mapsto\Lambda^+(\mu)$.
\normalsize
\end{center}

%%%%%%%%%%%%%%%%%%% END FIGURE %%%%%%%%%%%%%%%%%%%%%%%%%%%%%%%%%%%%%%%%%%%%%%%%

\medskip
As a consequence of Theorem~\ref{expmom}(ii), the maximum on the right-hand side of
\eqref{LegTra} is attained at some $\mu>-\rho(a^{**})$ if $b>b^{**}$ and at
$\mu=-\rho(a^{**})$ if $0\leq b\leq b^{**}$. Analogous assertions hold for
$\Lambda^-$, in particular, the restriction of $ I$ to $(-\infty,0]$ is the
Legendre transform of $\Lambda^-$. Since $\Lambda^-(\mu)=\Lambda^+(-\mu)$, the
moment generating function equals
    \begin{equation}
    \Lambda(\mu) = \lim_{T\to\infty}\frac 1T\log E\bigl(e^{- H_T}e^{\mu B_T}\bigr)
    =\max\{\Lambda^+(\mu),\Lambda^-(\mu)\}
    =\Lambda^+(|\mu|), \qquad \mu\in\R,
    \end{equation}
which is symmetric and strictly convex on $\R$, and non-differentiable at $0$, with
$\Lambda(0)=-a^*$ and $\lim_{\mu\da 0} \Lambda'(\mu)=b^*$.

\medskip
The outline of the present paper is as follows. In Section~\ref{sec-LDEM} we introduce
some preparatory material that will be needed in the sequel. Two basic propositions
are presented in Section~\ref{sec-basic}: a representation for the probabilities of
certain events under the Edwards measure, and an integrable majorant under which the
dominated convergence theorem can be applied. In Section~\ref{sec-proofsmainresults}
we carry out the proofs of Theorems~\ref{thm-LDEM}--\ref{expmom}. Some more refined
results about the Edwards model (which will be needed in a forthcoming paper \cite{HHK3})
appear in Section~\ref{sec-refine}. Finally, Section~\ref{sec-proof2} contains a
technical proof of a result used in Section~\ref{sec-refine}.

%%%%% SECTION 2 %%%%%%%%%%%%%%%%%%%%%%%%%%%%%%%%%%%%%%%%%

\section{Preliminaries}
\label{sec-LDEM}

In this section we provide some tools that are needed for the proofs of our main
results in Section~\ref{sec-EMtheorems}. These tools are taken from \cite{HH},
\cite{HHK1} and references cited therein. Section~\ref{analyt} introduces the
Sturm-Liouville operators that determine the constants. Section~\ref{BESQ}
provides the ingredients that are needed for the formulation of the Ray-Knight
Theorems (describing the joint distribution of the endpoint and the
local times), and contains a mixing property. Section~\ref{Airy} gives a spectral
decomposition of a function describing the ``overshoots'' of the path (i.e., the
pieces outside the interval between the starting point and the endpoint) in terms
of shifts of the Airy function, which plays an important role in various estimates.

\subsection{Sturm-Liouville operators and definition of the constants}
\label{analyt}

In \cite[Section~0.4]{HH} we introduced and analyzed a family of Sturm-Liouville
operators $\Kcal^a\colon L^2[0,\infty)\cap C^2[0,\infty)\to C[0,\infty)$,
indexed by $a\in\R$, defined as
    \begin{equation}
    \label{opdef}
    \bigl(\Kcal^a x\bigr)(h)=2hx''(h)+2x'(h)+(ah-h^2)x(h),\qquad h\geq 0.
    \end{equation}
The operator $\Kcal^a$ is symmetric and has a largest eigenvalue $\rho(a)\in\R$
with multiplicity one. The corresponding strictly positive (and $L^2$-normalized)
eigenfunction $x_a\colon [0,\infty)\to(0,\infty)$ is real-analytic and vanishes
faster than exponential at infinity, more precisely,
     \begin{equation}
     \label{xasyrough}
     \lim_{h\to\infty}h^{-{\frac 32}}\log x_a(h)=-\frac{\sqrt 2}3.
     \end{equation}

The eigenvalue function $\rho\colon\R\to\R$ has the following properties:
    \begin{equation}
    \label{rhodef}
    \begin{array}{ll}
    &\mbox{(a) } \rho \mbox{ is real-analytic};\\
    &\mbox{(b) } \rho \mbox{ is strictly log-convex, strictly convex and
    strictly increasing};\\
    &\mbox{(c) } \lim_{a\da-\infty}\rho(a)=-\infty, \rho(0)<0,
    \lim_{a\to\infty} \rho(a)=\infty.
    \end{array}
    \end{equation}
In terms of this object, the numbers $a^*,b^*,c^*$ appearing in
Theorem~\ref{thm-LDEM*} are defined as
    \begin{equation}
    \label{constdefs}
    \rho(a^*)=0,\qquad b^*=\frac 1{\rho'(a^*)},\qquad
    c^{*2}=\frac{\rho''(a^*)}{\rho'(a^*)^3},
    \end{equation}
while the numbers $a^{**},b^{**}$ appearing in Theorem~\ref{thm-LDEM} are
defined as
    \begin{equation}
    \label{b**def}
    a^{**}=2^{\frac 13}(-a_0), \qquad b^{**}=\frac 1{\rho'(a^{**})},
    \end{equation}
where $a_0$ ($\approx -2.3381$) is the largest zero of the Airy function:
    \begin{equation}
    \label{Airydiff}
    \begin{array}{ll}
    &\mbox{Ai is the unique solution of the Airy differential equation}\\
    &y''(h)=hy(h) \mbox{ that vanishes at infinity}.
    \end{array}
    \end{equation}
>From \cite[Lemma 6]{HHK1} we know that $a^*<-a_0$. Therefore $a^{**}>a^*$,
which in turn implies that $b^{**}<b^*$.

\subsection{Squared Bessel processes, a Girsanov transformation, and a
mixing property}
\label{BESQ}

The basic tools in our study of the Edwards model are the Ray-Knight
Theorems, which give a description of the joint distribution of the local
time process $(L(T,x))_{x\in\R}$ and the endpoint $B_T$. These will be
summarized in Proposition~\ref{repr} below. The key objects entering
into this description are introduced here.

The first key ingredients are:
    \begin{enumerate}
    \item[(i)]
    a squared two-dimensional Bessel process (BESQ$^2$),
    $X=(X_v)_{v\geq 0}$,
    \item[(ii)]
    a squared zero-dimensional Bessel process (BESQ$^0$),
    $X^\star=(X^{\star}_v)_{v\geq 0}$,
    \end{enumerate}
and their additive functionals
    \begin{equation}
    \label{Adef}
    A(t)=\int_0^t X_v\,\d v,\qquad A^\star(t)=\int_0^t X_v^{\star}\,\d v,\qquad t\geq 0.
    \end{equation}
The respective generators of BESQ$^2$ and BESQ$^0$ are given by
\footnote{BESQ$^0$ is sometimes called {\em Feller's diffusion}.}
     \begin{equation}
     Gf(h)=2hf''(h)+2f'(h), \qquad G^\star f(h)=2hf''(h),
     \end{equation}
for sufficiently smooth functions $f\colon[0,\infty)\to\R$. For $h\geq 0$, we write
$\P_h$ and $\P_h^\star$ to denote the probability law of $X$ and $X^\star$ given $X_0=h$
and $X^\star_0=h$, respectively. BESQ$^2$ takes values in $C^+=C^+[0,\infty)$, the
set of non-negative continuous functions on $[0,\infty)$. It has 0 as an entrance
boundary, which is not visited in finite positive time with probability one. BESQ$^0$
takes values in $C^+_0=C^+_0[0,\infty)$, the subset of those functions in $C^+$ that hit
zero and afterwards stay at zero. It has 0 as an absorbing boundary, which is visited
in finite time with probability one.

The second key ingredient is a certain Girsanov transformation, which turns
BESQ$^2$ into a diffusion with strong recurrence properties. Namely, the process
$(D_y^{\smallsup{a}})_{y\geq 0}$ defined by
    \begin{equation}
    \label{Ddef}
    D_y^{\smallsup{a}}=\frac{x_a(X_y)}{x_a(X_0)} \exp\Bigl\{-\int_0^y
    [(X_v)^2-aX_v+\rho(a)]\,\d v\Bigl\},\qquad y\geq 0,
    \end{equation}
is a martingale under $\P_h$ for any $h\geq 0$ and hence serves as a density
with respect to a new Markov process in the sense of a Girsanov transformation.
More precisely, the transformed process, which we also denote by $X=(X_v)_{v\geq 0}$,
has the transition density
    \begin{equation}
    \label{transdens}
    \widehat P^a_y(h_1,h_2)\,\d h_2
    =\E_{h_1}\bigl(D_y^{\smallsup{a}} \1_{\{X_y\in\d h_2\}}\bigr),\qquad y,h_1,h_2\geq 0.
    \end{equation}
We write $\widehat\P^a_h$ to denote the probability law of the transformed process
$X$ given $X_0=h$. This transformed process possesses the invariant distribution
$x_a(h)^2\,\d h$, and so
    \begin{equation}
    \label{widehatPdef}
    \widehat\P^a=\int_0^\infty \d h\,x_a(h)^2~\P^a_h
    \end{equation}
is its probability law in equilibrium. The transformed process is reversible under
$\widehat\P^a$, since BESQ$^2$ is reversible with respect to the Lebesgue
measure. Hence, $x_a(h_1)^2\widehat P^a_y(h_1,h_2)$ is symmetric in $h_1,h_2\geq 0$
for any $y\ge 0$.

The third key ingredient is the time-changed transformed process
$Y=X\circ A^{-1}=(X_{A^{-1}(t)})_{t\geq 0}$. We write $\widetilde\P^a_h$ to denote
the probability law of $Y$ given $Y_0=h$. This process
possesses the invariant distribution $\frac{1}{\rho'(a)}hx_a(h)^2\d h$, and so
    \begin{equation}
    \label{widetildePdef}
    \widetilde\P^a=\frac{1}{\rho'(a)} \int_0^\infty \d h~hx_a(h)^2~\widetilde\P^a_h
    \end{equation}
is its probability law in equilibrium. Both transformed processes
$X$ and $Y=X\circ A^{-1}$ are ergodic.

The following mixing property will be used frequently in the sequel. By
$\langle\cdot,\cdot\rangle$ we denote the inner product on $L^2=L^2[0,\infty)$,
and we write $\langle f,g\rangle_\circ =\int_0^\infty\d h\,h f(h)g(h)$ for the
inner product on $L^2$ weighted with the identity. The latter space will be
denoted by $L^{2,\circ}=L^{2,\circ}[0,\infty)$.

     \medskip
     \begin{prop}
     \label{convergence}
     Fix $a\in\R$ and fix measurable functions $f,g\colon[0,\infty)\to\R$ such
     that $f/{\rm id},g \in L^{2,\circ}$. For any family of measurable functions
     $f_s, g_s\colon[0,\infty)\to\R$, $s\geq 0$, such that $f_s/{\rm id},g_s \in
     L^{2,\circ}$, $s\geq 0$, and $f_s\to f, g_s\to g$ as $s\to\infty$ uniformly
     on compacts and in $L^{2,\circ}$, and for any family $a_s$, $s\geq 0$, such
     that $a_s\to a$ as $s\to\infty$,
     \begin{equation}
     \lim_{s\to\infty}\widehat\E^{a_s}\Bigl(\frac{f_s(X_0)}{x_{a_s}(X_0)}
     \frac{g_s(Y_{s})}{x_{a_s}(Y_{s})}\Bigr)
     =\langle f,x_a\rangle\frac 1{\rho'(a)}\langle g,x_a\rangle_\circ.
     \end{equation}
     \end{prop}

\noindent
This proposition is a slight extension of Proposition 3 in \cite{HHK1}; we omit the proof.

\subsection{BESQ$\boldsymbol{^0}$, the Airy function, and a spectral
decomposition}
\label{Airy}

For $a<a^{**}$, introduce the function $y_a\colon [0,\infty)\to(0,\infty]$
defined by
    \begin{equation}
    \label{ydef}
    y_a(h)=\E_h^\star\Bigl(e^{\int_0^\infty [aX_v^\star-(X_v^\star)^2]\,\d v}\Bigr).
    \end{equation}
(As a consequence of \eqref{wdef} and Proposition~\ref{lem-wasy} below,
the expectation on the right-hand side is infinite for $a>a^{**}$.)
It is known (see \cite{HHK1}, Lemma 5) that $y_a$ is equal to a normalized
scaled shift of the Airy function ${\rm Ai}$:
    \begin{equation}
        \label{yident}
        y_a(h)=\frac{{\rm Ai}(2^{-\frac13}(h-a))}{{\rm Ai}(-2^{-\frac13}a)}.
        \qquad h\geq 0,
        \end{equation}
It is well-known (see \cite[p.~43]{Erde56} and \eqref{Airyasy} below)
that $y_a$ vanishes faster than exponential at infinity:
     \begin{equation}
     \label{yasyrough}
     \lim_{h\to\infty}h^{-{\frac 32}}\log y_a(h)=-\frac{\sqrt 2}3.
     \end{equation}

An important role is played in the sequel by the function
$w\colon[0,\infty)^2\to[0,\infty)$ defined by
    \begin{equation}
    \label{wdef}
    w(h,t)\,\d t=\E_h^\star\Bigl(e^{-\int_0^\infty (X_v^\star)^2\,\d v}
    \1_{\{A^\star(\infty)\in \d t\}}\Bigr).
    \end{equation}
It is easily seen from \eqref{Adef} and \eqref{ydef} that $\int_0^\infty \d t\,
e^{at}w(h,t)=y_a(h)$ for $a<a^{**}$. We also have the following representation
for $w(h,t)$ derived in \cite[Lemma 7]{HHK1}:
\begin{eqnarray}
    \label{wident}
    w(h,t)&=&E_{\frac h2}\Bigl(e^{-2\int_0^t B_s\,\d s}\,\Big|\,T_0=t\Bigr)
    ~\varphi_h(t),\\
    \varphi_h(t)&=&\frac{P_{\frac h2}(T_0\in\d t)}{\d t}
    =(8\pi)^{-\frac12}t^{-\frac 32}h e^{-\frac{h^2}{8t}},\label{pdf}
    \end{eqnarray}
with $T_0=\inf\{t>0\colon B_t=0\}$ the first time $B$ hits zero. (We write $P_h$
and $E_h$ for probability and expectation with respect to standard Brownian
motion $B$ starting at $h\geq 0$, so that $P=P_0$, $E=E_0$.)

We will need the following expansion of the function $w$ in terms of shifts of
the Airy function:

    \medskip
    \begin{prop}
    \label{lem-wasy}
    \mbox{}
    \begin{itemize}
    \item[(i)]
    For any $\eps>0$,
        \begin{equation}
        \label{expansion}
        w(h,t)=\sum_{k=0}^\infty e^{a^{\smallsup{k}}(t-\eps)}\langle w(\cdot,\eps),
        \ef_k(\cdot)\rangle \,\ef_k(h),\qquad h\geq 0,\, t\geq\eps,
        \end{equation}
    where
        \begin{equation}
        \label{efdef}
        a^{\smallsup{k}} = 2^{\frac 13} a_k, \qquad
        \ef_k(h) = c_k\Ai(2^{-\frac13}(h+a^{\smallsup{k}})),\qquad h\geq 0,
        \end{equation}
    with $a_k$ the $k$-th largest zero of $\Ai$ and with $c_k$ chosen such
    that $\|\ef_k\|_2=1$.
    \item[(ii)]
    There exist constants $K_1,K_2,K_3\in(0,\infty)$ such that
        \begin{eqnarray}
        -a^{\smallsup{k}} &\sim& K_1 k^{\frac 23}, \qquad k\to\infty,
        \label{es1}\\
        \int_0^\infty h\ef_k(h)^2 \,\d h &\leq&  K_2 k^{\frac 23}\qquad  \forall\,k,
        \label{es2}\\
        \int_0^\infty \frac 1h \ef_k(h)^2 \,\d h &\leq & K_3 k^{\frac 13} \qquad \forall\,k.
        \label{es3}
        \end{eqnarray}
    \end{itemize}
    (Note that $a^{\smallsup{0}}=-a^{**}$ by (\ref{b**def}).)
    \end{prop}

\begin{proofsect}{Proof.}
(i) The proof comes in steps. We write $c$ for a generic constant in $(0,\infty)$
whose value may change from appearance to appearance.

\medskip
{\bf 1.} Let $\Kcal^\star$ be the second-order differential operator on $C_0^\infty=
C_0^\infty[0,\infty)$, the set of smooth functions $x\colon [0,\infty)\to\R$
that vanish at zero, defined by
    \begin{equation}
    \label{Kstardef}
    (\Kcal^\star x)(h)=\begin{cases}2x''(h)-hx(h)&\mbox{if }h>0,\\
    0&\mbox{if }h=0.
    \end{cases}
    \end{equation}
This operator is symmetric with respect to the $L^2$-inner product on $L^2_0 = L^2\cap
C_0^\infty$. Furthermore, we can identify all the eigenvalues and eigenfunctions
of $\Kcal^\star$ in $L^2_0$ in terms of scaled shifts of the Airy function. Namely,
a comparison of \eqref{Airydiff} and \eqref{Kstardef} shows that the $k$-th eigenspace
is spanned by the eigenfunction $\ef_k\colon [0,\infty)\to\R$ given in \eqref{efdef}
and the $k$-th eigenvalue is $a^{\smallsup{k}}$, $k\in\N_0$.

\medskip
{\bf 2.} We next show that $\Kcal^\star$ has a compact inverse on $L^2$. Therefore,
this inverse has an orthonormal basis of eigenvectors in $L^2$, and hence the same
is true for $\Kcal^\star$ itself. Consequently, $(\ef_k)_{k\in\N_0}$ is an orthonormal
basis of $L^2$. This fact will be needed later.

We begin by identifying the inverse of $\Kcal^\star$. To do so, we follow
\cite{G81}. Let
    \begin{equation}
    \label{ydefs}
    y_1(u)= \Bi(2^{\frac 13}u)-\Bi(0)\frac{\Ai(2^{\frac 13}u)}{\Ai(0)},
    \qquad y_2(u)=\Ai(2^{\frac 13}u),
    \end{equation}
where $\Ai$ is the Airy function and $\Bi$ is another, linearly independent,
solution to \eqref{Airydiff} (for the precise definitions of $\Ai$ and $\Bi$,
see \cite[10.4.1--10.4.3]{AS70}). Hence, both $y_1$ and $y_2$ solve
$\Kcal^\star y=0$, $y_1$ satisfies the boundary condition at zero ($y_1(0)=0$),
while $y_2$ satisfies the boundary condition at infinity ($y_2\in L^2$). Let
$G\colon[0,\infty)^2\to \R$ (Green function) be defined by
    \begin{equation}
    \label{zdef}
    G(u,v) = K y_1(u\wedge v) y_2(u \vee v) \quad \mbox{ with }
    \quad K=-2y_1'(0)y_2(0).
    \end{equation}
Let $\Gamma$ be the operator on $L^2$ defined by
    \begin{equation}
    (\Gamma y)(u) = \int_0^{\infty} G(u,v) y(v) \,\d v.
    \end{equation}
According to \cite[Proposition 2.15]{G81}, $x=\Gamma y$ is a weak solution of
the equation $\Kcal^\star x=y$ with boundary condition $x(0)=0$, for any $y\in L^2$.
In fact, we can adapt the proof of \cite[Proposition 9.12]{G81} to see that $\Gamma$ is
the inverse of $\Kcal^\star$, since $\Kcal^\star x=0$ does not have solutions
in $L^2$ that satisfy the boundary condition $x(0)=0$. Hence, we are done once
we show that $\Gamma$ is a compact operator.

\medskip
{\bf 3.} By \cite[Theorem 8.54]{G81}, it suffices to show that $\Gamma$ is a
Hilbert-Schmidt operator, i.e., $G$ is square-integrable on $[0,\infty)^2$. In
order to show this, we first note that \eqref{zdef} gives
    \begin{equation}
    \label{HS2}
    \int_0^\infty \d u \int_0^\infty \d v\,\, G^2(u,v)
    = 2 K^2 \int_0^\infty \d u \int_0^u \d v\,\, y_2(u)^2 y_1(v)^2.
    \end{equation}
Substitute \eqref{ydefs} to see that, since $\Ai\in L^2$, it suffices to show that
    \begin{equation}
    \int_0^\infty \d u \int_0^u \d v\,\,
    \Ai(u)^2\Bi(v)^2 < \infty.
    \end{equation}
Since $\Bi$ is locally bounded and $\Ai\in L^2$, the latter amounts to
    \begin{equation}
    \int_1^\infty \d u \int_1^u \d v\,\,
    \Ai(u)^2\Bi(v)^2 < \infty.
    \end{equation}
We next use \cite[10.4.59 and 10.4.63]{AS70}, which shows that
    \begin{equation}
    \Ai(u) \leq c u^{-\frac14} e^{-\frac 23 u^{\frac 32}}, \qquad
    \Bi(v) \leq c v^{-\frac14} e^{\frac 23 v^{\frac 32}},\qquad u,v\geq 1.
    \label{asyAiBi}
    \end{equation}
Hence
    \begin{equation}
    \int_1^\infty \d u \int_1^u \d v\,\, \Ai(u)^2\Bi(v)^2
    \leq c^4 \int_1^\infty \d u\, u^{-\frac 12} \int_1^u \d v\, v^{-\frac 12}
    e^{-\frac 43 (u^{\frac 32}-v^{\frac 32})}.
    \end{equation}
Use partial integration to see that
    \begin{equation}
    \int_1^u \d v\, v^{-\frac 12} e^{-\frac 43 (u^{\frac 32}-v^{\frac 32})}
    = \frac 12 \int_1^u \d v\, v^{-1} \frac{\d}{\d v}\Bigl(
    e^{-\frac 43 (u^{\frac 32}-v^{\frac 32})}\Bigr)\leq
    \frac 12 \big[v^{-1}e^{-\frac 43 (u^{\frac 32}-v^{\frac 32})}\big]_{v=1}^u
    \leq \frac 12 u^{-1}, \qquad u\geq 1.
    \end{equation}
Hence
    \begin{equation}
    \int_1^\infty \d u \int_1^u \d v\,\,\Ai(u)^2\Bi(v)^2
    \leq \frac 12 c^4 \int_1^\infty \d u\, u^{-\frac 32} < \infty.
    \end{equation}
This proves that $\Gamma$ is a compact operator, so that $(\ef_k)_{k\in \N_0}$
is an orthonormal basis of $L^2$.

\medskip
{\bf 4.} To prove the expansion in \eqref{expansion}, we now need the following:

    \medskip
    \begin{lemma}
    \label{bvp}
    For any $\eps>0$, the function $w$ is a solution of the initial-boundary-value
    problem
    \begin{equation}
    \label{bvp*}
    \begin{array}{rcll}
    \partial_t w(h,t)&=&\Kcal^\star \bigl(w(\cdot,t)\bigr)(h),
    \qquad &h\ge 0, t>\eps,\\
    w(0,t)&\equiv& 0,\qquad &t\ge\eps,
    \end{array}
    \end{equation}
    and the initial value $w(\cdot,\eps)$ lies in $C_0^\infty$.
    \end{lemma}

\begin{proofsect}{Proof.}
Use the Markov property at time $s>0$ in \eqref{wident} to see that, for any
$h>0$ and $t>s$,
    \begin{equation}
    w(h,t)=E_{\frac h2}\Bigl(e^{-\int_0^s 2B_v\,\d v}\1_{\{T_0>s\}}w(2B_s,t-s)\Bigr).
    \end{equation}
Now differentiate with respect to $s$ at $s=0$, to obtain
    \begin{equation}
    0=-h w(h,t)+2(\partial_h)^2 w(h,t)-\partial_t w(h,t)
    =\Kcal^\star \bigl(w(\cdot,t)\bigr)(h)-\partial_t w(h,t).
    \end{equation}
This shows that the partial differential equation in \eqref{bvp*} is satisfied on
$(0,\infty)^2$. It is clear that it is also satisfied at the boundary where $h=0$,
since $w(0,t)=0$ for all $t>0$ (recall (\ref{wdef}--\ref{wident})).
\end{proofsect}
\qed

\medskip
{\bf 5.} From \eqref{wident} it follows that $w(\cdot,\eps)\in C_0^\infty$ for any
$\eps>0$. A spectral decomposition in terms of the eigenvalues
$(a^{\smallsup{k}})_{k\in\N_0}$ and the eigenfunctions $(\ef_k)_{k\in\N_0}$ of
$\Kcal^\star$ shows that \eqref{bvp*} has the solution given in \eqref{expansion}.

\medskip\noindent
(ii) In \cite[10.4.94,10.4.96,10.4.97,10.4.105]{AS70} the following asymptotics
for the Airy function can be found. As $k\to\infty$,
\begin{equation}
\label{Airyasymp}
- a_k \sim c k^{\frac 23}, \qquad a_{k-1}-a_k \sim c k^{-\frac 13},
\qquad \max_{[a_k,a_{k-1}]} |\Ai| \sim ck^{-\frac 16},
\qquad |\Ai'(a_k)| \sim c k^{\frac 16}.
\end{equation}
We will use these in combination with the observation that, by \eqref{Airydiff},
$\Ai$ is convex (concave) between any two successive zeroes where it is negative
(positive).

The first assertion in \eqref{Airyasymp} is \eqref{es1}. To prove
(\ref{es2}--\ref{es3}), we write the recursion
    \begin{equation}
    \label{csplit}
    c_k^{-2}=\int_{0}^{\infty} \Ai\bigl(2^{-\frac13}(h+a^{\smallsup{k}})\bigr)^2\,\d h
    =c_{k-1}^{-2}+2^{\frac 13}\int_{a_k}^{a_{k-1}}\Ai(h)^2\,\d h.
    \end{equation}
Using the second and third assertion in \eqref{Airyasymp}, we find that
$\int_{a_k}^{a_{k-1}}\Ai(h)^2\,\d h \asymp k^{-\frac 23}$ and hence that
$c_k^{-2}\asymp k^{\frac 13}$. In a similar way, we find that
    \begin{equation}
    \label{split}
    \int_{0}^{\infty} h \Ai\bigl(2^{-\frac13}(h+a^{\smallsup{k}})\bigr)^2\,\d h
    \leq c k, \qquad
    \int_{0}^{\infty} \frac 1h \Ai\bigl(2^{-\frac13}(h+a^{\smallsup{k}})\bigr)^2\,\d h
    \leq c k^{\frac 23}.
\end{equation}
Combining \eqref{split} with \eqref{efdef} and $c_k^{-2}\asymp k^{\frac 13}$,
we obtain (\ref{es2}--\ref{es3}).
\end{proofsect}
\qed

%%%%%%%%%%%% SECTION 3 %%%%%%%%%%%%%%%%

\section{Two basic propositions}
\label{sec-basic}

In this section we present the basic tools of our proofs. Section \ref{RKdescr}
introduces the Ray-Knight Theorems, which give a flexible representation for the
probabilities of certain events under the Edwards measure. Section \ref{domi}
exhibits an integrable majorant under which limits may be interchanged with
integrals.

\subsection{Ray-Knight representation}
\label{RKdescr}

In this section we formulate the Ray-Knight Theorems that were already
announced in Section~\ref{BESQ}. We do this in the compact form derived in
\cite[Section 1.2]{HHK1}, which is best suited for the arguments in the sequel.

For any measurable set $G\subset C_0$, define $w_G\colon [0,\infty)\times
[0,\infty)\to\R$ by
    \begin{equation}
    \label{Wdef}
    w_G(h,t)\,\d t=\E_h^\star\Bigl(e^{-\int_0^\infty (X_v^\star)^2\,\d v}
    \1_{\{X^\star\in G\}}\1_{\{A^\star(\infty)\in \d t\}}\Bigr).
    \end{equation}
It is clear that $w_G$ is increasing in $G$. For $G=C_0$, $w_{C_0}$ is identical
to $w$ defined in \eqref{wdef}.

For $y\geq0$, denote by $C^+[0,y]$ the set of non-negative continuous functions
on $[0,y]$. Then the set $\Ccal^+=\bigcup_{y\geq 0}\bigl(\{y\}\times C^+[0,y]\bigr)$
is the appropriate state space of the pair $\bigl(B_T,L(T,B_T-\cdot)|_{[0,B_T]}\bigr)$
consisting of the endpoint $B_T(\geq 0)$ and the local time process between the
endpoint $B_T$ and the starting point 0.

    \medskip
    \begin{prop}[Ray-Knight representation]
    \label{repr}
    Fix $a\in\R$. Then, for any $T>0$ and any measurable sets $G^+,G^-\subset C_0$ and
    $F\subset  \Ccal$,
    \begin{equation}
    \label{repr1}
    \begin{aligned}
    e^{aT} E\Bigl(e^{- H_T}&e^{-\rho(a)B_T}\1_{\{ L(T,B_T+\cdot\,)\in G^+\}}
    \1_{\{(B_T,L(T,B_T-\cdot\,)|_{[0,B_T]})\in F\}}\1_{\{L(T,-~\cdot\,)\in G^-\}}\Bigr)\\
    &=\int_0^\infty\d t_1\int_0^\infty\d t_2\, \1_{\{t_1+t_2\leq T\}}e^{a(t_1+t_2)}\\
    &\quad\times\widehat\E^a\Bigl(\1_{\{(A^{-1}(T-t_1-t_2),X|_{[0,A^{-1}(T-t_1-t_2)]})\in F\}}
    \frac{w_{G^+}(X_0,t_1)}{x_a(X_0)}
    \frac{w_{G^-}(Y_{T-t_1-t_2},t_2)}{x_a(Y_{T-t_1-t_2})}\Bigr).
    \end{aligned}
    \end{equation}
    \end{prop}

\begin{proofsect}{Proof.}
We briefly indicate how \eqref{repr1} comes about. Details can be found in
\cite[Section 1.2]{HHK1}. Recall the notation in Section~\ref{BESQ}. Fix $T>0$.
Then, according to the Ray-Knight Theorems, for any $t_1,t_2,h_1,h_2\geq 0$
and $y>0$, conditioned on the event
    \begin{equation}
    \label{event}
    \{B_T=y\}\cap\{L(T,B_T)=h_1\}\cap\{L(T,0)=h_2\}
    \cap\Bigl\{\int_{B_T}^\infty L(T,x)\,\d x=t_1\Bigr\}
    \cap\Bigl\{\int_0^\infty L(T,-x)\,\d x=t_2\Bigr\},
    \end{equation}
the joint distribution of the processes
    \begin{equation}
    L(T,B_T+\cdot\,), \qquad L(B_T-\,\cdot\,)|_{[0,y]}, \qquad L(T,-\,\cdot\,),
    \end{equation}
on $C^+_0\times C^+[0,y]\times C^+_0$ is equal to the joint distribution of the processes
    \begin{equation}
    X^{\star,1}(\cdot), \qquad X(\cdot)|_{[0,y]}, \qquad X^{\star,2}(\cdot),
    \end{equation}
under
    \begin{equation}
    \label{Measure}
    \P^\star_{h_1}(\,\cdot\,|A^\star(\infty)=t_1)
    \otimes\P_{h_1}(\,\cdot\,|A(y)=T-t_1-t_2,X_y=h_2)
    \otimes\P^\star_{h_2}(\,\cdot\,|A^\star(\infty)=t_2),
    \end{equation}
where $X$ is BESQ$^2$ and $X^{\star,1}$, $X^{\star,2}$ are independent copies of
BESQ$^0$. In particular, the intersection local time in \eqref{HT} has the
representation
    \begin{equation}
    \label{HTident}
     H_T\stackrel{\rm{law}}{=}\int_0^\infty
    (X_v^{\star,1})^2\,\d v+\int_0^y (X_v)^2\,\d v
    + \int_0^\infty (X_v^{\star,2})^2\,\d v.
    \end{equation}
Use \eqref{transdens} for $y=A^{-1}(T-t_1-t_2)$ and note that, on the event
$\{A(T-t_1-t_2)=y\}\cap\{X_0=h_1,X_y=h_2\}$, \eqref{Ddef} becomes
    \begin{equation}
    D_y^{\smallsup{a}}=\frac{x_a(h_2)}{x_a(h_1)}
    \exp\Bigl\{-\int_0^y (X_v)^2\,\d v\Bigr\}
    e^{a(T-t_1-t_2)}e^{-\rho(a)y},
    \end{equation}
which implies that
    \begin{equation}
    e^{aT}e^{- H_T} e^{-\rho(a)B_T}\stackrel{\rm{law}}
    {=}\frac{x_a(h_1)}{x_a(h_2)}D_y^{\smallsup{a}}e^{a(t_1+t_2)} e^{-\int_0^\infty
    (X_v^{\star,1})^2\,\d v}e^{-\int_0^\infty
    (X_v^{\star,2})^2\,\d v}.
    \end{equation}
Integrate the left-hand side with respect to $P$ and the right-hand side with respect
to the measure in \eqref{Measure}, and absorb the term $D_y^{\smallsup{a}}$ into
the notation of the transformed diffusion. Integrate over $h_1,h_2\geq 0$ and note
that $X_0$ has the distribution $x_a(h_1)^2\,\d h_1$ under $\widehat \E^a$. Finally,
use the notation in \eqref{Wdef}, to obtain \eqref{repr1}.
\end{proofsect}
\qed

\subsection{Domination}
\label{domi}

In order to perform the limit $T\to\infty$ on the right-hand side of \eqref{repr1},
we will need the dominated convergence theorem to interchange this limit
with the integrals over $t_1$ and $t_2$. The following proposition provides the
required domination.

    \medskip
    \begin{prop}[Domination]
    \label{dominate}
    For any $a_s$, $s\geq 0$, in a compact subset of $(-\infty,a^{**})$, the map
    \begin{equation}
    \label{domibound}
    (t_1,t_2)\mapsto \sup_{s\geq 0}~ e^{a_s(t_1+t_2)}
    \widehat\E^{a_s}\Bigl(\frac{w(X_0,t_1)}{x_{a_s}(X_0)}
    \frac{w(Y_{s},t_2)}{x_{a_s}(Y_{s})}\Bigr)
    \end{equation}
    is integrable over $(0,\infty)^2$.
    \end{prop}

\begin{proofsect}{Proof.}
Under the expectation in \eqref{domibound} we make a change of measure from the
invariant distribution of $X$ to the invariant distribution of $Y$, i.e., we
replace $\widehat\E^{a_s}$ by $\widetilde\E^{a_s}$ and add a factor of $\rho'(a_s)/Y_0$.
Fix $1<p\leq q<\infty$ such that $\frac {1}{p}+\frac {1}{q}=1$, apply H\"older's
inequality and use the stationarity of $Y$ under $\widetilde\P^{a_s}$. This gives,
for any $t_1,t_2>0$, the bound
    \begin{equation}
    \label{Hoelder}
    \widetilde\E^{a_s}\Bigl(\frac{w(Y_0,t_1)}{Y_0 x_{a_s}(Y_0)}
    \frac{w(Y_{s},t_2)}{x_{a_s}(Y_{s})}\Bigr)\leq
    W^{\smallsup{1}}_{p}(t_1)W^{\smallsup{2}}_{q}(t_2),
    \end{equation}
where  the functions $W^{\smallsup{1}}_{p},W^{\smallsup{2}}_{q}\colon(0,\infty)
\to(0,\infty)$ are defined by
    \begin{equation}
    \label{Wdefs}
    W^{\smallsup{1}}_{p}(t)=\widetilde\E^{a_s}\Big(
    \Big(\frac{w(Y_0,t)}{Y_0x_{a_s}(Y_0)}\Big)^{p}\Big)^{\frac 1p},\qquad
    W^{\smallsup{2}}_{q}(t)=\widetilde\E^{a_s}\Big(
    \Big(\frac{w(Y_0,t)}{x_{a_s}(Y_0)}\Big)^{q}\Big)^{\frac 1q}.
    \end{equation}
Hence, it suffices to show that the maps
    \begin{equation}
    t\mapsto e^{a_st} W^{\smallsup{1}}_{p}(t),
    \qquad t\mapsto e^{a_st}W^{\smallsup{2}}_{q}(t),
    \end{equation}
are integrable at zero and at infinity, uniformly in $s$, for a suitable choice
of $p$ and $q$. In the proof of Proposition~4 in \cite{HHK1} we showed that
$W^{\smallsup{1}}_{p}$ and $W^{\smallsup{2}}_{q}$, with $a_s$ replaced by $a^*$,
are integrable at zero when
$p<q$ with $p,q$ sufficiently close to 2. An inspection of the proof shows that
they are actually integrable at zero uniformly in $s$.

We will show that $t\mapsto e^{a_s t}
W^{\smallsup{1}}_{2}(t)$ and $t\mapsto e^{a_s t}W^{\smallsup{2}}_{2}(t)$ are
integrable at infinity uniformly in $s$. This will complete the proof because the left-hand
side of \eqref{Hoelder} does not depend on $p,q$.

We use Proposition~\ref{lem-wasy} with $\eps=1$ together with the
representations (recall \eqref{widetildePdef})
    \begin{equation}
    W^{\smallsup{1}}_{2}(t)
    =\frac 1{\sqrt{\rho'(a_s)}} \Bigl(\int_0^\infty \d h\, \frac 1h
    w(h,t)^{2}\Bigr)^{\frac{1}{2}}, \qquad
    W^{\smallsup{2}}_{2}(t)
    = \frac 1{\sqrt{\rho'(a_s)}}\Bigl(\int_0^\infty \d h\, h
    w(h,t)^{2}\Bigr)^{\frac{1}{2}}.
    \end{equation}
Using \eqref{expansion}, the Cauchy-Schwarz inequality and the fact that $\|\ef_k\|_2=1$,
we estimate
\begin{equation}
    W^{\smallsup{1}}_{2}(t)\leq\frac{1}{\sqrt{\rho'(a_s)}}
    \Bigl(\|w(\cdot,1)\|_2^2
    \sum_{k_1, k_2=0}^\infty e^{(a^{\smallsup{k_1}}+a^{\smallsup{k_2}})(t-1)}
    \int_0^\infty \frac 1h |\ef_{k_1}(h)| |\ef_{k_2}(h)|~\d h\,\Bigr)^{\frac{1}{2}},
    \qquad t\geq 1.
    \end{equation}
Using the Cauchy-Schwarz inequality for the last integral, we obtain the bound
    \begin{equation}\label{W1esti}
    W^{\smallsup{1}}_{2}(t)\leq\frac{\|w(\cdot,1)\|_2}{\sqrt{\rho'(a_s)}}
    \sum_{k=0}^\infty e^{a^{\smallsup{k}}(t-1)}
    \Big(\int_0^\infty \frac 1h \ef_{k}(h)^2~\d h\Big)^{\frac 12},\qquad t\geq 1.
    \end{equation}
In the same way, we find that
    \begin{equation}\label{W2esti}
    W^{\smallsup{2}}_{2}(t)
    \leq\frac{\|w(\cdot,1)\|_2}{\sqrt{\rho'(a_s)}}
    \sum_{k=0}^\infty e^{a^{\smallsup{k}}(t-1)}
    \Big(\int_0^\infty h \ef_{k}(h)^2~\d h\Big)^{\frac 12},\qquad t\geq 1.
    \end{equation}
Substitute (\ref{es2}--\ref{es3}) into (\ref{W1esti}--\ref{W2esti}) and use that
$a^{\smallsup{k}}\leq a^{\smallsup{0}}=-a^{**}$, to estimate
    \begin{equation}
    W^{\smallsup{1}}_{2}(t)\vee W^{\smallsup{2}}_{2}(t)
    \leq c\, e^{-a^{**}(t-2)}\sum_{k=0}^\infty e^{a^{\smallsup{k}}}k^{\frac 13},
    \qquad t\geq 2.
    \end{equation}
By \eqref{es1}, the sum in the right-hand side converges. Since $a_s<a^{**}$, $s\geq 0$,
is bounded away from $a^{**}$, it is now obvious that the maps $t\mapsto e^{a_st}
W^{\smallsup{1}}_{2}(t)$ and $t\mapsto e^{a_st}W^{\smallsup{2}}_{2}(t)$ are integrable
at infinity uniformly in $s$.
\end{proofsect}
\qed

%%%%%%%%%%%%%%%%%%%%%%%%%%%% SECTION 4 %%%%%%%%%%%%%%%%%%%%%%%%%%%%%%%%%%%

\section{Proof of Theorems~\ref{thm-LDEM}--\ref{expmom}}
\label{sec-proofsmainresults}

In Sections \ref{expmomi-ii}--\ref{expmomi-ii*} we give the proof of
Theorems~\ref{thm-LDEM}--\ref{expmom} with the help of Propositions
\ref{repr}--\ref{dominate}. In Section \ref{sec-growth} we derive a
technical proposition that is needed along the way.

\subsection{Growth rate of a restricted moment generating function}
\label{sec-growth}

Abbreviate $B_{[0,T]}=\{B_t\colon t\in[0,T]\}$ for the range of the path up to
time $T$. For $T>0$ and $\delta,C\in(0,\infty]$, define events
    \begin{eqnarray}
    \eve(\delta;T)
    &=& \bigl\{B_{[0,T]}\subset[-\delta,B_T+\delta]\bigr\},
    \label{Ehatdef}\\
    \evf^{\leq}(\delta,C;T)
    &=&\Big\{\max_{x\in[-\delta,\delta]}
    L(T,x)\leq C,\max_{x\in[B_T-\delta,B_T+\delta]}L(T,x)\leq C\Big\}.
    \label{Fhatdeflq}
\end{eqnarray}
In words, on $\eve(\delta;T)$ the path does not visit more than the $\delta$-neighborhood
of the interval between its starting point 0 and its endpoint $B_T$, while on
$\evf^{\leq}(\delta,C;T)$ its local times in the $\delta$-neighborhoods
of these two points are bounded by $C$. Note that both $\eve(\infty;T)$ and
$\evf^{\leq}(\delta,\infty;T)$ are the full space.

    \medskip
    \begin{prop}
    \label{mainparti}
    Fix $\mu>-\rho(a^{**})$. Then, for any $\delta,C\in(0,\infty]$
    there exists a constant $K_1(\delta,C)\in(0,\infty)$ such that,
    for any $\mu_T\to\mu$ as $T\to\infty$,
    \begin{equation}\label{refine1}
    e^{\rho^{-1}(-\mu_T)T} E\Bigl(e^{- H_T}e^{\mu_T B_T}\1_{\eve(\delta,T)}
    \1_{\evf^\leq(\delta,C;T)}\1_{\{B_T\geq 0\}}\Bigr) = K_1(\delta,C)+o(1).
    \end{equation}
    Moreover, if $\mu=\mu_b$ solves $I(b)=\mu b-\Lambda^+(\mu)$, then
    the same is true when $\1_{\{B_T\geq 0\}}$ is replaced by
    $\1_{\{B_T\approx bT\}}$.
    \end{prop}

\begin{proofsect}{Proof.}
We may assume that $\mu_T>-\rho(a^{**})$ for all $T$. Fix $\delta,C\in(0,\infty]$
and choose $a_T$ such that $\mu_T+\rho(a_T)=0$, i.e., $a_T=\rho^{-1}(-\mu_T)<a^{**}$.
Clearly, $\lim_{T\to\infty}a_T=\rho^{-1}(-\mu)<a^{**}$.
Since, on $\eve(\delta;T)\cap\{B_T\leq 2\delta\}$, we can estimate
\begin{equation}
H_T=4\delta\int_{-\delta}^{3\delta}\frac{\d x}{4\delta}\, L(T,x)^2
\geq 4\delta\Bigl(\int_{-\delta}^{3\delta}\frac{\d x}{4\delta}\, L(T,x)\Bigr)^2
=\frac {T^2}{4\delta},
\end{equation}
we may insert the indicator of $\{B_T\geq 2\delta\}$ in the expectation on
the left-hand side of \eqref{refine1}, paying only a factor $1+o(1)$ as $T\to\infty$.

{\bf 1.}
Introduce the following subsets of $C^+_0$, respectively, $\Ccal^+$ (see below \eqref{Wdef}):
    \begin{eqnarray}
    G^\leq_{\delta,C}&=&\bigl\{g\in C^+_0\colon g(\delta)=0,\max g\leq C\bigr\},
    \label{GdeltaCdef}\\
    F^\leq_{\delta,C}&=&\Bigl\{(y,f)\in\Ccal^+\colon y\geq 2\delta,
    \max_{[0,\delta]}f\leq C,
    \max_{[y-\delta,y]}f\leq C\Bigr\}.
    \label{FdeltaCdef}
    \end{eqnarray}
Note that
\begin{equation}
\begin{array}{ll}
&\eve(\delta;T) \cap \evf^\leq (\delta,C;T)\cap \{B_T\geq 2\delta\}\\[0.2cm]
&\qquad = \{L(T,B_T+\cdot) \in G^\leq_{\delta,C}\}
\cap \{L(T,-\,\cdot) \in G^\leq_{\delta,C}\}
\cap \{(B_T,L(T,B_T-\cdot)|_{[0,B_T]}) \in F^\leq_{\delta,C}\}.
\end{array}
\end{equation}
Apply Proposition~\ref{repr} for $a=a_T$ with $F=F^\leq_{\delta,C}$ and
$G^+=G^-=G^\leq_{\delta,C}$, to get
    \begin{equation}
    \label{refine1a}
    \begin{aligned}
    \mbox{l.h.s.~of \eqref{refine1}}
    &=(1+o(1))\int_0^\infty\d t_1\int_0^\infty\d t_2\, \1_{\{t_1+t_2\leq T\}}e^{a_T(t_1+t_2)}\\
    &\qquad\times \widehat\E^{a_T}\left(
    \frac{w_{G^\leq_{\delta,C}}(X_0,t_1)}{x_{a_T}(X_0)}
    \1_{\{A^{-1}(T-t_1-t_2)\geq 2\delta\}}
    \1_{\bigl\{\max_{[0,\delta]}X\leq C\bigr\}}\right.\\
    &\qquad\times
    \left. \1_{\bigl\{\max_{[A^{-1}(T-t_1-t_2)-\delta,A^{-1}(T-t_1-t_2)]}X\leq C\bigr\}}
    \frac{w_{G^\leq_{\delta,C}}(Y_{T-t_1-t_2},t_2)}{x_{a_T}(Y_{T-t_1-t_2})}\right).
    \end{aligned}
    \end{equation}

{\bf 2.}
In the case $C=\infty$, the last two indicators vanish and we can identify the limit
of the integrand as $T\to\infty$ with the help of Lemma~\ref{convergence}. Indeed,
apply Lemma~\ref{convergence} for $f(\cdot)=w_{G_{\delta}}(\cdot,t_1)$ and
$g(\cdot)=w_{G_{\delta}}(\cdot,t_2)$, where we put $G_\delta=G^\leq_{\delta,\infty}
=\{g\in C^+_0\colon g(\delta)=0\}$. Then we obtain that the integrand converges to
    \begin{equation}
    e^{a(t_1+t_2)}\langle w_{G_{\delta}}(\cdot,t_1),x_a\rangle\frac 1{\rho'(a)}
    \langle w_{G_{\delta}}(\cdot,t_2),x_a\rangle_\circ,
    \end{equation}
where we also use that $A^{-1}(\infty)=\infty$ because $X$ never hits 0 (recall
\eqref{Adef}). According to Proposition~\ref{dominate}, we are allowed to
interchange the limit $T\to\infty$ with the two integrals over $t_1$ and $t_2$.
This implies that \eqref{refine1} holds with $K_1(\delta,\infty)$ identified as
    \begin{equation}
    \label{K1ident}
    K_1(\delta,\infty)
    =\langle y_a^{\smallsup{\delta}},x_a\rangle
    \frac1{\rho'(a)}\langle y_a^{\smallsup{\delta}},x_a\rangle_\circ,
    \end{equation}
where $y_a^{\smallsup{\delta}}(h)$ is defined as (recall \eqref{Wdef})
    \begin{equation}
    \label{ydeltadef}
    y_a^{\smallsup{\delta}}(h)=\int_0^\infty\d t\,e^{at}w_{G_\delta}(h,t)
    =\E_h^\star\Bigl(e^{\int_0^\infty[aX_v^\star- (X_v^\star)^2]\,\d v}
    \1_{\{X^\star_\delta=0\}}\Bigr).
    \end{equation}
Trivially, $K_1(\delta,\infty)>0$. Since $y_a^{\smallsup{\delta}}\leq y_a$, it
follows from \eqref{xasyrough} and \eqref{yasyrough} that $K_1(\delta,\infty)<\infty$.

{\bf 3.}
Next we return to \eqref{refine1a} and consider the case $C\in(0,\infty)$.
Note that the integrals over $t_1$ and $t_2$ can both be restricted to
$[0,C\delta]$, since $w_{G^\leq_{\delta,C}}(h,t)=0$ for $t>C\delta$ as is
seen from \eqref{Wdef} and \eqref{GdeltaCdef}.

Let us abbreviate $s=T-t_1-t_2$. We first apply the Markov property for
the process $X$ at time $\delta$ and integrate over all values $z=A(\delta)$.
Because of the appearance of the indicator of $\{\max_{[0,\delta]}X\leq C\}$, we may restrict
to $z\in[0,C\delta]$ (recall \eqref{Adef}). We note that the additive functional
of the process $(X_{\delta+t})_{t\geq 0}$ given that $A(\delta)=z$, denoted
by $\widetilde A = (\widetilde A(t))_{t\geq 0}$, is given by $\widetilde
A(t)=A(t+\delta)-z$. Making the change of variables $s=\widetilde A(t)+z$,
we see that $A^{-1}(s)=\widetilde A^{-1}(s-z)+\delta$ for any $s\ge 0$. Defining
$f^{t_1}_{s,T}\colon(0,\infty)^2\to[0,\infty)$ by
    \begin{equation}
    \label{fTdef}
    \begin{aligned}
    %\mbox{}&
    \!\!\!\!f^{t_1}_{s,T}(h,z)\,\d h\d z%\\
    &=x_{a_T}(h)
    \widehat\E^{a_T}\left(\frac{w_{G^\leq_{\delta,C}}(X_0,t_1)}
    {x_{a_T}(X_0)}\1_{\{A^{-1}(s)\geq 2\delta\}}
    \1_{\bigl\{\max_{[0,\delta]}X\leq C\bigr\}}
    \1_{\{X_\delta\in \d h\}}\1_{\{A(\delta)\in \d z\}}\right),
   \end{aligned}
   \end{equation}
we thus obtain that the expectation under the integral in \eqref{refine1a}
can be written as
    \begin{equation}
    \label{refine1d}
    \begin{aligned}
    \widehat\E^{a_T}&\Bigl(
    \frac{w_{G^\leq_{\delta,C}}(X_0,t_1)}{x_{a_T}(X_0)}
    \1_{\{A^{-1}(s)\geq 2\delta\}}
    \1_{\bigl\{\max_{[0,\delta]}X\leq C\bigr\}}
    \1_{\bigl\{\max_{[A^{-1}(s)-\delta,A^{-1}(s)]}X\leq C\bigr\}}
    \frac{w_{G^\leq_{\delta,C}}(Y_{s},t_2)}{x_{a_T}(Y_{s})}\Bigr)\\
    &=\int_0^{C\delta} \d z\, \widehat \E^{a_T}\Bigl(\frac{f^{t_1}_{s,T}(X_0,z)}{x_{a_T}(X_0)}
    \1_{\bigl\{\max_{[A^{-1}(s-z),A^{-1}(s-z)+\delta]}X\leq C\bigr\}}
    \frac{w_{G^\leq_{\delta,C}}(X_{A^{-1}(s-z)+\delta},t_2)}{x_{a_T}
    (X_{A^{-1}(s-z)+\delta})}\Bigr).
    \end{aligned}
    \end{equation}
(The tilde can be removed afterwards.) We next apply the Markov property
for the process $Y$ at time $s-z$ (respectively, the strong Markov property for
the process $X$ at time $A^{-1}(s-z)$), to write
    \begin{equation}
    \label{refine1e}
    \mbox{r.h.s.~of \eqref{refine1d}}
    =\int_0^{C\delta} \d z\, \widehat \E^{a_T}\Bigl(\frac{f^{t_1}_{s,T}(X_0,z)}{x_{a_T}(X_0)}
    \frac{g^{t_2}_T(Y_{s-z})}{x_{a_T}(Y_{s-z})}\Bigr),
    \end{equation}
where $g^{t_2}_T$ is defined by
    \begin{equation}
    \label{gTdef}
    g^{t_2}_T(h)=x_{a_T}(h)\widehat \E^{a_T}_h\Bigl(\1_{\bigl\{\max_{[0,\delta]}X\leq C\bigr\}}
    \frac{w_{G^\leq_{\delta,C}}(X_{\delta},t_2)}{x_{a_T}(X_{\delta})}\Bigr).
    \end{equation}

{\bf 4.}
We want to take the limit $s\to\infty$ in \eqref{refine1e} (recall that $s=T-t_1-t_2$)
and use Proposition \ref{convergence}. Therefore we need dominated convergence. To
establish this, we note that
     \begin{equation}
     \label{supest}
     \sup_{h\in[0,C]}~\sup_{t\in[0,C\delta]}~\sup_{T\ge 1}~\frac{w(h,t)}{x_{a_T}(h)}
     = K < \infty
     \end{equation}
(see (\ref{wdef}--\ref{pdf}) and recall that $x_a$ is bounded away from zero on
$[0,C]$ and continuous in $a$). By (\ref{gTdef}--\ref{supest}), the last quotient
in the right-hand side of \eqref{refine1e} is bounded above by $K$. Substituting
\eqref{fTdef} into \eqref{refine1e} and using that $w_{G^\leq_{\delta,C}}
\leq w_{C^+_0} = w$, we therefore obtain
    \begin{equation}
    \label{refineextra2}
    \begin{aligned}
    \mbox{integrand of r.h.s.~of \eqref{refine1e}}
    &\le K\, \widehat\E^{a_T}\Bigl(
    \frac{f^{t_1}_{s,T}(X_0,z)}{x_{a_T}(X_0)}\Bigr)\\
    &\le K ~\frac{\widehat\E^{a_T}\Bigl(
    \frac{w(X_0,t_1)}{x_{a_T}(X_0)}
    \1_{\bigl\{\max_{[0,\delta]}X\leq C\bigr\}}
    \1_{\{A(\delta)\in \d z\}}\Bigr)}{\d z}\\
    & \le K^2~\frac{\widehat\P^{a_T}\bigl(A(\delta)\in\d z\bigr)}{\d z}.
    \end{aligned}
    \end{equation}
It is easy to see from \eqref{Ddef} that the right-hand side of \eqref{refineextra2}
is bounded uniformly in $T\geq 1$ and $z\in [0,C\delta]$. Therefore we have an
integrable majorant for \eqref{refine1e}, which allows us to interchange the limit
$s\to\infty$ with the integral over $z$.

{\bf 5.}
In order to identify the limit as $s\to\infty$ of the integrand on the right-hand
side of \eqref{refine1e}, we apply Lemma~\ref{convergence} to see that this
integrand converges to $\langle f^{t_1}(\cdot,z),x_a(\cdot)\rangle\frac 1{\rho'(a)}
\langle g^{t_2},x_a\rangle_\circ$, with $f^{t_1}$ and $g^{t_2}$ the pointwise limit
of $f^{t_1}_{s,T}$ and $g^{t_2}_T$, respectively:
    \begin{eqnarray}
    \label{fdef}
    f^{t_1}(h,z)\,\d h \d z &=& x_a(h)\widehat\E^a\Bigl(\frac{w_{G^\leq_{\delta,C}}(X_0,t_1)}
    {x_a(X_0)} \1_{\bigl\{\max\limits_{[0,\delta]}X\leq C\bigr\}}\1_{\{X_\delta\in \d h\}}
    \1_{\{A(\delta)\in\d z\}}\Bigr)\\
    g^{t_2}(h)\,\d h &=& x_a(h)\widehat\E^a_h\Bigl(\1_{\bigl\{\max_{[0,\delta]}X\leq C\bigr\}}
    \frac{w_{G^\leq_{\delta,C}}(X_0,t_2)}{x_a(X_\delta)}\Bigr).
    \end{eqnarray}
Using this in \eqref{refine1e} and interchanging the integral over $z$ with the limit
$s\to\infty$, we obtain that
    \begin{equation}
    \label{refine1g}
    \lim_{s\to\infty}\bigl(\mbox{l.h.s.~of \eqref{refine1d}}\bigr)
    =\langle f^{t_1},x_a\rangle\frac 1{\rho'(a)}\langle
    g^{t_2},x_a\rangle_\circ
    \end{equation}
with $f^{t_1}(h) = \int_0^{C\delta} \d z\,f^{t_1}(h,z)$.

{\bf 6.}
Finally, recall that $s=T-t_1-t_2$ and that $e^{a(t_1+t_2)}$ times the left-hand side
of \eqref{refine1d} is equal to the integrand on the right-hand side of \eqref{refine1a}.
According to Proposition~\ref{dominate}, we are allowed to interchange the limit
$T\to\infty$ with the two integrals over $t_1$ and $t_2$. Hence we obtain that
\eqref{refine1} holds with $K_1(\delta,C)$ identified as the integral over $t_1,t_2$
of the right-hand side of \eqref{refine1g}, which is a strictly positive finite number.
This proves the statement with the indicator on $\1_{\{B_T\geq 0\}}$.

{\bf 7.} To prove the statement with $\1_{\{B_T\geq 0\}}$ replaced by
$\1_{\{B_T\approx bT\}}$, we let $\mu=\mu_b$ solve $I(b)=\mu b-\Lambda^+(\mu)$.
The statement follows when we show that for every $a\in \R$, we have that
    \begin{equation}\label{aima}
    e^{\rho^{-1}(-\mu)T} E\Bigl(e^{a\frac{B_T-bT}{\sqrt{T}}}e^{- H_T}e^{\mu B_T}\1_{\eve(\delta,T)}
    \1_{\evf^\leq(\delta,C;T)}\1_{\{B_T\geq 0\}}\Bigr) = e^{\frac{a^2}{2} \sigma_b^2}
    K_1(\delta,C)+o(1)
    \end{equation}
for some $\sigma_b^2\in (0,\infty)$. Indeed, \eqref{aima} shows that
$\1_{\{|B_T-bT|>\gamma_T, B_T\geq 0\}}$ is asymptotically negligible
for any $\gamma_T$ such that $\gamma_T/\sqrt{T} \rightarrow \infty$.

In order to prove \eqref{aima}, we rewrite the left-hand side as
    \begin{equation}\label{aimb}
    e^{[\rho^{-1}(-\mu)-\rho^{-1}(-\mu_{a,T})]T- a b\sqrt{T}}
    e^{\rho^{-1}(-\mu_{a,T})T}E\Bigl(e^{-H_T}e^{\mu_{a,T} B_T}\1_{\eve(\delta,T)}
    \1_{\evf^\leq(\delta,C;T)}\1_{\{B_T\geq 0\}}\Bigr),
    \end{equation}
where $\mu_{a,T}=\mu+\frac{a}{\sqrt{T}}$. Clearly, $\mu_{a,T}\rightarrow \mu$, so that
the second factor converges to $K_1(\delta,C)$. We are therefore left to compute the
exponential. We note that since $\mu=\mu_b$ solves $I(b)=\mu b-\Lambda^+(\mu)$, we have
that $\rho'(-\mu_b)=1/b$. Therefore,
    \begin{equation}
    \rho^{-1}(-\mu_{a,T})=\rho^{-1}(-\mu)-\frac{a}{\sqrt{T}} \frac{1}{\rho'(-\mu)}+
    \frac{a^2}{2T} \frac{\d^2 }{\d \mu^2} \rho^{-1}(-\mu) +o(T^{-1}).
    \end{equation}
Therefore,
    \begin{equation}
    e^{[\rho^{-1}(-\mu)-\rho^{-1}(-\mu_{a,T})]T- a b\sqrt{T}} =
    e^{\frac{a^2}{2} \frac{\d^2 }{\d \mu^2} \rho^{-1}(-\mu)}(1+o(1)),
    \end{equation}
which completes the proof with $\sigma_b^2=-\frac{\d^2 }{\d \mu^2} \rho^{-1}(-\mu_b)$.
    \end{proofsect}
\qed

\subsection{Proof of Theorem \ref{expmom}(i--iii)}
\label{expmomi-ii}

\setcounter{step}0

    \medskip
    \begin{step}
    \label{lem-Lambda+ident}
    For any $\mu>-\rho(a^{**})$, the limit in \eqref{Lambda+def} exists and equals
    $\Lambda^+(\mu)=-\rho^{-1}(-\mu)$. On $(-\rho(a^{**}),\infty)$, the function
    $\Lambda^+$ is real-analytic and strictly convex, and satisfies
    $\lim_{\mu\da -\rho(a^{**})}(\Lambda^+)'(\mu)=b^{**}$.
    \end{step}

\begin{proofsect}{Proof.}
Fix $\mu>-\rho(a^{**})$, apply Proposition~\ref{mainparti} with $\delta=C=\infty$, and use
the continuity of $\rho$, to obtain that the limit in the definition of $\Lambda ^+(\mu)$
in \eqref{Lambda+def} exists and equals $-\rho^{-1}(-\mu)$. This proves the first assertion.
The remaining assertions follow from (\ref{rhodef}--\ref{b**def}).
\end{proofsect}
\qed

In the following step, we consider paths that never go below $-\delta$, have local times
that are bounded by $C$ in the $\delta$-neighborhood of the starting point 0, and have
the endpoint $B_T$ close to 0. Recall that $\gamma_T$ is a
function that satisfies $\gamma_T/T\to0$ and $\gamma_T/\sqrt T\to\infty$ as $T\to\infty$.

    \medskip
    \begin{step}
    \label{zero}
    For any $\delta\in(0,\infty)$ and $C\in(0,\infty]$,
    \begin{equation}
    \label{Izero}
     E\Bigl(e^{- H_T}\1_{\{B_T\in[0,\gamma_T]\}}
    \1_{\bigl\{\min_{[0,T]}B\geq -\delta\bigr\}}
    \1_{\bigl\{\max_{[-\delta,\delta]}L(T,\cdot)\leq C\bigr\}}\Bigr)
    \geq e^{-a^{**}T+o(T)}, \qquad T\to\infty.
    \end{equation}
    \es

\begin{proofsect}{Proof.}
Pick $a=a^{**}$ and apply Proposition~\ref{repr} for
    \begin{equation}
    F=F_{\delta,C}=\{(y,f)\in\Ccal^+\colon y\leq \delta,\max_{[y-\delta,y]}f\leq C\},
    \qquad
    G^+=C^+_0,
    \qquad
    G^-=G^\leq_{\delta,C}
    \end{equation}
(recall \eqref{GdeltaCdef}). Note that the
event under the expectation on the left-hand side of \eqref{Izero} contains the event
\begin{equation}
\{L(T,B_T+\cdot)\in C^+_0\} \cap \{(B_T,L(T,B_T-\cdot))\in F_{\delta,C}\}
\cap \{L(T,-\cdot)\in G^\leq_{\delta,C}\}.
\end{equation}
Also note that $e^{-\rho(a^{**})B_T}\leq 1$ when $B_T\geq 0$ because $\rho(a^{**})>0$.
Therefore we find
    \begin{equation}
    \begin{aligned}
    \mbox{l.h.s.~of \eqref{Izero}}
    &\geq \int_0^\infty \d t_1 \int_0^\infty \d t_2\,\1_{\{t_1+t_2\leq T\}}e^{-a^{**}s}\\
    &\qquad\times\widehat\E^{a^{**}}\Bigl(\frac{w(X_0,t_1)}{x_{a^{**}}(X_0)}
    \1_{\{A^{-1}(s)\leq\delta\}}\1_{\big\{\max_{[A^{-1}(s)-\delta,A^{-1}(s)]}X\leq C\big\}}
    \frac{w_{G^\leq_{\delta,C}}(Y_{s},t_2)}{x_{a^{**}}(Y_{s})}\Bigr),
    \end{aligned}
    \end{equation}
where we again abbreviate $s=T-t_1-t_2$. Next we interchange the two integrals,
restrict the $t_2$-integral to $[0,\delta]$ and the $t_1$-integral to
$[T-t_2-\delta,T-t_2]$, estimate $A^{-1}(s)\leq A^{-1}(\delta)$ for $s\leq\delta$,
and integrate over $s=T-t_1-t_2$, to get
    \begin{equation}
    \begin{aligned}
    \mbox{}&\mbox{l.h.s.~of \eqref{Izero}}\\
    &\geq \int_0^\delta \d t_2\int_0^\delta \d s\,
    \widehat\E^{a^{**}}\Bigl(\frac{w(X_0,T-t_2-s)}{x_{a^{**}}(X_0)}
    \1_{\{A^{-1}(\delta)\leq\delta\}}
    \1_{\{\max_{[0,\delta]}X\leq C\}}
    \frac{w_{G^\leq_{\delta,C}}(Y_{s},t_2)}{x_{a^{**}}(Y_{s})}\Bigr).
    \end{aligned}
    \end{equation}
Now we use Proposition~\ref{lem-wasy}(i) to estimate $w(X_0,T-s-t_2)\geq e^{-a^{**}T+o(T)}$,
uniformly on the domain of integration. The remaining expectation on the right-hand side
no longer depends on $T$ and is strictly positive for any $\delta\in(0,\infty)$ and
$C\in(0,\infty]$.
\end{proofsect}
\qed

     \medskip
     \bes
     \label{Lambdaconst}
     $\Lambda^+$ equals $-a^{**}$ on $(-\infty,-\rho(a^{**})]$.
     \es

\begin{proofsect}{Proof.}
For $\mu\leq -\rho(a^{**})$, define $\Lambda^+_-(\mu)$ and $\Lambda^+_+(\mu)$ as in
\eqref{Lambda+def} with $\lim$ replaced by $\liminf$ and $\limsup$, respectively.
Since $\Lambda^+_+$ is obviously non-decreasing, we have $\Lambda^+_+(\mu)\leq
\Lambda^+(-\rho(a^{**})+\eps)$ for $\mu\leq -\rho(a^{**})$ and any $\eps>0$.
Using Step~\ref{lem-Lambda+ident} and the continuity of $\rho$, we see that
$\lim_{\eps\da 0}\Lambda^+(-\rho(a^{**})+\eps)=-\rho^{-1}(\rho(a^{**}))=-a^{**}$,
which shows that $\Lambda^+_+(\mu)\leq -a^{**}$. In order to get the reversed
inequality for $\Lambda^+_-(\mu)$, bound
    \begin{equation}
     E\Bigl(e^{- H_T}e^{\mu B_T}\1_{\{B_T\geq0\}}\Bigr)\geq E\Bigl(e^{- H_T}e^{\mu B_T}
    \1_{\{B_T\in[0, \gamma_T]\}}\Bigr)
    \geq e^{\mu\gamma_T} E\Bigl(e^{- H_T}\1_{\{B_T\in[0, \gamma_T]\}}\Bigr),
    \end{equation}
take logs, divide by $T$, let $T\to\infty$ and use Step~\ref{zero}, to obtain
that $\Lambda^+_-(\mu)\geq -a^{**}$. Since $\Lambda^+_-\leq \Lambda^+_+$, this
implies the assertion.
\end{proofsect}
\qed

    \medskip
    \begin{step}
    \label{Lambda+lim}
    $\Lambda^+(\mu)=\frac 12 \mu^2+\Ocal(\mu^{-1})$ as $\mu\to\infty$.
    \es

\begin{proofsect}{Proof.}
According to Step~\ref{lem-Lambda+ident}, we have $\Lambda^+(\mu)=-\rho^{-1}(-\mu)$
for $\mu>-\rho(a^{**})$. Hence, in order to obtain the asymptotics for
$\Lambda^+(\mu)$ as $\mu\to\infty$, we need to obtain the asymptotics for $\rho(a)$
as $a\to -\infty$. In the following we consider $a<0$.

We use Rayleigh's Principle (see \cite[Proposition 10.10]{G81}) to write
(recall \eqref{opdef})
    \begin{equation}
    \begin{array}{lll}
    \rho(a) &=& \sup\limits_{x\in L^2\cap C^2 \colon \|x\|_2=1}
    \langle \Kcal^a x,x\rangle\\[0.2cm]
    &=& \sup\limits_{x\in L^2\cap C^2 \colon \|x\|_2=1}
    \int_0^{\infty} \big[-2hx'(h)^2 +(ah-h^2) x(h)^2\big]~\d h.
    \end{array}
    \end{equation}
Substituting $x(h)=(-a)^{\frac14} y((-a)^{\frac 12} h)$, we get
    \begin{equation}
    \label{rhoscal}
    \rho(a) = (-a)^{\frac 12} \sup_{y\in L^2\cap C^2 \colon \|y\|_2=1}
    \int_0^{\infty} \big[-2hy'(h)^2 -(h+h^2(-a)^{-\frac 32}) y(h)^2\big]~\d h.
    \end{equation}
Hence, we have the upper bound $\rho(a)\leq V(-a)^{\frac 12}$ with
    \begin{equation}\label{Vdef}
    V =\sup_{y\in L^2\cap C^2 \colon \|y\|_2=1}
    \int_0^{\infty}\big[-2hy'(h)^2 -hy(h)^2\big]~\d h.
    \end{equation}
By completing the square under the integral and partially integrating the cross
term, we easily see that $y^*(h)= \frac{1}{\sqrt{2}}e^{-h/\sqrt{2}}$ is the maximizer
of \eqref{Vdef} and $V=-\sqrt 2$. Substituting $y^*$ into \eqref{rhoscal}, we can
also bound $\rho(a)$ from below:
    \begin{equation}
    \rho(a)
    \geq  -\sqrt{2}\,(-a)^{\frac 12}-(-a)^{-1}\int_{0}^{\infty} h^2 y^*(h)^2~\d h.
    \end{equation}
Therefore,
    \begin{equation}
    \rho(a) = -\sqrt{2}\,(-a)^{\frac 12} + \Ocal(|a|^{-1}),\qquad a\to-\infty.
    \end{equation}
Consequently,
    \begin{equation}
    \Lambda^+(\mu)=-\rho^{-1}(-\mu)=\frac{1}{2}\mu^2 +\Ocal(\mu^{-1}),
     \qquad\mu\to\infty.
    \end{equation}
\end{proofsect}
\qed

Steps \ref{lem-Lambda+ident}, \ref{Lambdaconst} and \ref{Lambda+lim} complete the proof of
Theorem~\ref{expmom}(i--iii).

\subsection{Proof of Theorem~\ref{thm-LDEM} and \ref{expmom}(iv)}
\label{expmomi-ii*}

For $b\in\R$, define $ I_-(b)$ and $ I_+(b)$ as in \eqref{rateED} with
$\lim$ replaced by $\limsup$ and $\liminf$, respectively.

    \medskip
    \bes
    \label{LegTraproof}
    For any $b>b^{**}$, the limit in \eqref{rateED} exists and \eqref{LegTra} holds.
    \es

\begin{proofsect}{Proof.}
Fix $b>b^{**}$.

{\bf 1.}
To derive \lq$\geq$\rq\ in \eqref{LegTra} for $ I_-$ instead of $ I$, bound, for
any $\mu\in\R$,
    \begin{equation}
    \label{widehatIlb}
    \begin{aligned}
     E\bigl(e^{-H_T} \1_{\{|B_T-b T|\leq \gamma_T\}}\bigr)
    &\leq e^{-\mu bT+|\mu|\gamma_T} E\bigl(e^{-H_T}e^{\mu B_T}
    \1_{\{|B_T-b T|\leq \gamma_T\}}\bigr)\\
    &\leq e^{-\mu bT+|\mu|\gamma_T}  E\bigl(e^{-H_T}e^{\mu B_T}
    \1_{\{B_T\geq 0\}}\bigr),
    \end{aligned}
    \end{equation}
where the last inequality holds for any $T$ sufficiently large
because $\gamma_T/T\to 0$ as $T\to\infty$. Take logs, divide by $T$, let
$T\to\infty$, use \eqref{Lambda+def} and minimize over $\mu\in\R$, to obtain
    \begin{equation}
    \label{Iminbound}
    - I_-(b) \leq \min_{\mu\in\R} [-\mu b+\Lambda^+(\mu)].
    \end{equation}
This shows that \lq$\geq$\rq\ holds in \eqref{LegTra} for $I$ replaced by
$I_-$.

{\bf 2.}
To derive \lq$\leq$\rq~in \eqref{LegTra}, bound, for any $\mu\in\R$,
    \begin{equation}
    \label{trickproof1}
    \begin{aligned}
     E\bigl(&e^{- H_T} \1_{\{|B_T-b T|\leq \gamma_T\}}\bigr)\\
    &\geq  E\bigl(e^{- H_T}\1_{\eve(\delta,T)}
    \1_{\{|B_T-b T|\leq \gamma_T\}}\1_{\{B_T\geq 0\}}\bigr)\\
    & \geq e^{-\mu bT-|\mu|\gamma_T}
     P^{\mu,\delta,T}\bigl(|B_T-b T|\leq \gamma_T\bigr)
     E\bigl(e^{- H_T}e^{\mu B_T}\1_{\eve(\delta,T)}
    \1_{\{B_T\geq 0\}}\bigr),
    \end{aligned}
\end{equation}
where $ P^{\mu,\delta,T}$ denotes the probability law whose density with respect
to $ P$ is proportional to $e^{- H_T}e^{\mu B_T}\1_{\eve(\delta,T)}\1_{\{B_T\geq 0\}}$.

{\bf 3.}
Let $\mu_b$ be the maximizer of the map $\mu\mapsto \mu b-\Lambda^+(\mu)$. (Note that,
by Step~\ref{lem-Lambda+ident}, the maximizer is unique and is characterized by
$(\Lambda^+)'(\mu_b)=b$.) Next we argue that
    \begin{equation}
    \label{trick}
    \lim_{T\to\infty}  P^{\mu_b,\delta,T}\bigl(|B_T-b T|\leq \gamma_T\bigr)=1.
    \end{equation}
Indeed, pick $\eps_T=\gamma_T/cT>0$ (with $c>0$ to be specified later) and estimate
    \begin{equation}
    \1_{\{B_T\geq bT+\gamma_T\}}\leq e^{\eps_T[B_T-b T-\gamma_T]}.
    \end{equation}
This implies, with the help of Step \ref{lem-Lambda+ident} and Proposition~\ref{mainparti}
with $\mu_T=\mu+\epsilon_T$, $C=\infty$, that
    \begin{equation}
    \label{trickproof}
     P^{\mu_b,\delta,T}\bigl(B_T\geq b T+ \gamma_T\bigr)
    \leq e^{-\eps_T[bT+\gamma_T]}e^{[\Lambda^+(\mu_b+\eps_T)
    -\Lambda^+(\mu_b)]T}(1+o(1)),\qquad T\to\infty.
    \end{equation}
A Taylor expansion of $\Lambda^+$ around $\mu_b$, in combination with the observation that
$(\Lambda^+)'(\mu_b)=b$ and $c=(\Lambda^+)''(\mu_b)>0$, yields that the right-hand
side of \eqref{trickproof} is equal to
    \begin{equation}
    e^{\frac c2 \eps_T^2T[1+\Ocal(\eps_T)]-\eps_T\gamma_T}
    =e^{-\frac{\gamma_T^2}{2cT}[1+\Ocal(\frac{\gamma_T}{T})]},\qquad T\to\infty.
    \end{equation}
The right-hand side vanishes as $T\to\infty$ because $\gamma_T/T\to 0$ and
$\gamma_T/\sqrt T\to\infty$. This shows that $\lim_{T\to\infty} P^{\mu_b,\delta,T}
\bigl(B_T\geq b T+ \gamma_T\bigr)=0$. Analogously, replacing $\eps_T$ by $-\eps_T$,
we can prove that $\lim_{T\to\infty} P^{\mu_b,\delta,T}\bigl(B_T\leq bT-\gamma_T\bigr)=0$.
Hence, \eqref{trick} holds.

{\bf 4.}
Use \eqref{trick} in \eqref{trickproof1} for $\mu=\mu_b$, take logs, divide
by $T$, let $T\to\infty$, and use Step~\ref{lem-Lambda+ident} and Proposition~\ref{mainparti},
to obtain
    \begin{equation}
    \label{Iplusbound}
    - I_+(b) \geq -\mu_b b + \Lambda^+(\mu_b)
    = - \max_{\mu\in\R}[\mu b-\Lambda^+(\mu)].
    \end{equation}
This shows that \lq$\leq$\rq\ holds in \eqref{LegTra} for $I$ replaced by
$I_+$. Combine \eqref{Iminbound} and \eqref{Iplusbound} to obtain that $ I_-=I=I_+$ and
that \eqref{LegTra} holds on $(b^{**},\infty)$.
\end{proofsect}
\qed

   \medskip
   \bes\label{<b**>}
   For any $b\geq 0$, $ I_-(b)\geq -b\rho(a^{**})+a^{**}$.
   \es

\begin{proofsect}{Proof.}
Estimate
\begin{equation}
\1_{\{|B_T-bT|\leq\gamma_T\}} \leq \1_{\{B_T \leq bT+\gamma_T\}} \leq
e^{-\rho(a^{**})[B_T-bT-\gamma_T]},
\end{equation}
to obtain, for $T$ sufficiently large,
    \begin{equation}
    \begin{aligned}
    E\Bigl(e^{-H_T}\1_{\{|B_T-bT|\leq\gamma_T\}}\Bigr)
    &\leq 2 E\Bigl(e^{-H_T}\1_{\{|B_T-bT|\leq\gamma_T\}}\1_{\{B_T\geq 0\}}\Bigr)\\
    &\leq 2 e^{b\rho(a^{**})T+\gamma_T \rho(a^{**})} E\Bigl(e^{-H_T}
    e^{-\rho(a^{**}) B_T}\1_{\{B_T\geq 0\}}\Bigr).
    \end{aligned}
    \end{equation}
According to the definition of $\Lambda^+$ in \eqref{Lambda+def}, the expectation
in the right-hand side is equal to $e^{\Lambda^+(-\rho(a^{**}))T+o(T)}$. We
therefore obtain that $I(b)\geq -b\rho(a^{**})-\Lambda^+(-\rho(a^{**}))$.
Now Step~\ref{Lambdaconst} concludes the proof.
\end{proofsect}
\qed

   \medskip
   \bes\label{<b**<}
   For any $0\leq b\leq b^{**}$, $I_+(b)\leq -b\rho(a^{**})+a^{**}$.
   \es

\begin{proofsect}{Proof.}
Fix $0\leq b\leq b^{**}$, pick $b'>b^{**}$ and put $\alpha=b/b'\in[0,1)$.
We split the path $(B_s)_{s\in[0,T]}$ into two pieces: $s\in[0,\alpha T]$
and $s\in[\alpha T,T]$. First we bound from below by inserting several indicators:
    \begin{equation}
    \begin{aligned}
     E\Bigl(&e^{- H_T}\1_{\{|B_T-bT|\leq \gamma_T\}}\Bigr)\\
    &\geq  E\Bigl(e^{- H_T}\1_{\{|B_{\alpha T}-b'\alpha T|\leq \gamma_T/2\}}
    \1_{\bigl\{\max_{[0,\alpha T]}B\leq B_{\alpha T}+\delta\bigr\}}
    \1_{\bigr\{\max_{[B_{\alpha T}-\delta,B_{\alpha T}+\delta]}
    L(\alpha T,\cdot)\leq C\bigr\}}\\
    &\qquad \times
    \1_{\{|\widetilde B_{(1-\alpha)T}|\leq \gamma_T/2\}}
    \1_{\bigl\{\min_{[0,(1-\alpha )T]}\widetilde B\geq -\delta]\bigr\}}
    \1_{\bigl\{\max_{[B_{\alpha T}-\delta,B_{\alpha T}+\delta]}
    \widetilde L((1-\alpha)T,\cdot)\leq C\bigr\}}\Bigr).
    \end{aligned}
    \end{equation}
Here, $(\widetilde B_s)_{s\in [0,(1-\alpha)T]}$ is the Brownian motion
with $\widetilde B_s = B_{\alpha T+s}-B_{\alpha T}$, and
$\widetilde L((1-\alpha)T,x) = L(T,x)-L(\alpha T,x), x\in \R,$
are its local times.

On the event under the expectation in the right-hand side, we may estimate
    \begin{equation}
    H_T= H_{\alpha T}+ \widetilde H_{(1-\alpha)T}
    +2\int_{B_{\alpha T}-\delta}^{B_{\alpha T}+\delta}
    L(\alpha T,x)\widetilde L((1-\alpha)T,x)\,\d x
    \leq  H_{\alpha T}+ \widetilde H_{(1-\alpha)T}+4\delta C^2,
    \end{equation}
where $ \widetilde H_{(1-\alpha)T}$ denotes the intersection local time for
the second piece. Using the Markov property at time $\alpha T$, we therefore
obtain the estimate
    \begin{equation}
    \begin{aligned}
    \mbox{}& E\Bigl(e^{- H_T}\1_{\{|B_T-bT|\leq \gamma_T\}}\Bigr)\\
    &\geq e^{-4\delta C^2} E\Bigl(e^{- H_{\alpha T}}
    \1_{\{|B_{\alpha T}-b'\alpha T|\leq \gamma_T/2\}}
    \1_{\bigl\{\max_{[0,\alpha T]}B\leq B_{\alpha T}+\delta\bigr\}}
    \1_{\bigl\{\max_{[B_{\alpha T}-\delta,B_{\alpha T}+\delta]}
    L(\alpha T,\cdot)\leq C\bigr\}}\Bigr)\\
    &\times  E\Bigl(e^{- H_{(1-\alpha) T}}
    \1_{\bigl\{|B_{(1-\alpha )T}|\leq \gamma_T/2\bigr\}}
    \1_{\bigl\{\min_{[0,(1-\alpha)T]}B\geq -\delta\bigr\}}
     \1_{\bigl\{\max_{[-\delta,\delta]} L((1-\alpha )T,\cdot)\leq C\bigr\}}\Bigr).
    \end{aligned}
    \end{equation}
(The tilde can be removed afterwards.) Now use Proposition~\ref{mainparti} (in
combination with an argument like in parts 2-3 of the proof of Step~\ref{LegTraproof})
for the first term (with $T$ replaced by $\alpha T$) and use Step~\ref{zero}
for the second term (with $T$ replaced by $(1-\alpha)T$), to conclude that
    \begin{equation}
    I(b)\leq \alpha  I(b')+\bigl(1-\alpha\bigr) a^{**}
    =\frac b{b'}\bigl( I(b')-a^{**}\bigr)+a^{**}.
    \end{equation}
Let $b'\da b^{**}$, use the continuity of $ I$ in $b^{**}$, and note that
$I(b^{**})-a^{**}=-b^{**}\rho(a^{**})$ by Step \ref{LegTraproof}, to conclude
the proof.
\end{proofsect}
\qed

    \medskip
    \bes
    Theorems~\ref{thm-LDEM} and \ref{expmom}(iv) hold.
    \es

\begin{proofsect}{Proof.}
Steps~\ref{lem-Lambda+ident} and \ref{LegTraproof} allow us to identify
$I$ on $(b^{**},\infty)$ as $I(b)= -b\rho(a_b)+a_b$, where $a_b$ solves
$\rho'(a_b)=1/b$ (the maximum in \eqref{LegTra} is attained at $\mu=-\rho(a_b)$).
{}From this and (\ref{rhodef}--\ref{b**def}) it follows that
    \begin{equation}
    I'(b) = -\rho(a_b), \qquad I''(b) = -\rho'(a_b)\frac{\d}{\d b}a_b
    = \frac{[\rho'(a_b)]^3}{\rho''(a_b)}>0, \qquad b>b^{**}.
   \end{equation}
In particular, $ I$ is real-analytic and strictly convex on $(b^{**},\infty)$.
Since $a_{b^{**}}=a^{**}$, it in turn follows that
    \begin{equation}
    \min_{b\geq 0} I(b)
    =\min_{b>b^{**}} I(b) = I(b^*) = a^*,
    \end{equation}
where $a^*$ solves $\rho(a^*)=0$ (the minimum is attained at $b^*=1/\rho'(a^*)$).
This, together with Steps \ref{LegTraproof}--\ref{<b**<}, proves
Theorem~\ref{thm-LDEM}(i--iii).

Step~\ref{LegTraproof} shows that \eqref{LegTra} holds on $(b^{**},\infty)$.
To show that it also holds on $[0,b^{**}]$, use Step~\ref{Lambdaconst} to get
    \begin{equation}
    -b\rho(a^{**})+a^{**} =\max_{\mu\in\R}[b\mu-\Lambda^+(\mu)],
    \qquad 0\leq b\leq b^{**},
    \end{equation}
since the maximum is attained at $\mu=-\rho(a^{**})$. Recall from Steps~\ref{<b**>}--\ref{<b**<} that the left-hand side is equal to $I(b)$. Thus we have proved
Theorem~\ref{expmom}(iv).

Finally, Theorem~\ref{thm-LDEM}(iv) is an immediate consequence of
Theorem~\ref{expmom}(iii--iv).
\end{proofsect}
\qed

%%%%%%%%%%%%%% SECTION 5 %%%%%%%%%%%%%%%%%%%%%%%%%%%%%%%%%%%%%%%%

\section{Addendum 1: An extension of Proposition~\ref{mainparti}}
\label{sec-refine}

At this point we have completed the proof of the main results in Section~\ref{EDW}.
In Sections~\ref{sec-refine}--\ref{sec-proof2} we derive an extension of
Proposition~\ref{mainparti} that will be needed in a forthcoming paper \cite{HHK3}.
In that paper we show that several one-dimensional polymers models in discrete
space and time, such as the weakly self-avoiding walk, converge to the Edwards
model, after appropriate scaling, in the limit of vanishing self-repellence.
The proof is based on a coarse-graining argument, for which we need
Proposition~\ref{mainpart} below.

Recall the events in (\ref{Ehatdef}--\ref{Fhatdeflq}). For
$\delta\in (0,\infty), \alpha\in[0,\infty)$,
define the event
    \begin{equation}
    \label{Fhatdefgq}
    \evf^{\,\geq}(\delta,\alpha;T)
    =\Big\{\max_{x\in[B_T-\delta,B_T+\delta]}
    L(T,x)\geq\alpha\delta^{-\frac 12}\Big\}.
   \end{equation}
Note that $\evf^{\,\geq}(\delta,0;T)$ is the full space.

    \medskip
    \begin{prop}
    \label{mainpart}
    Fix $\mu>-\rho(a^{**})$. Then:
    \begin{enumerate}
    \item[(i)]
    For any $\delta\in(0,\infty)$ and $\alpha\in[0,\infty)$ there exists a
    $K_2(\delta,\alpha)\in(0,\infty)$ such that
    \begin{equation}\label{refine2}
    e^{\rho^{-1}(-\mu)T} E\Bigl(e^{- H_T}e^{\mu B_T}\1_{\eve(\delta,T)}
    \1_{\evf^{\,\geq}(\delta,\alpha;T)}\1_{\{B_T\geq 0\}}\Bigr)
    = K_2(\delta,\alpha)+o(1), \qquad T\to\infty.
    \end{equation}
    \item[(ii)] For any $\alpha\in(0,\infty)$,
    \begin{equation}\label{Kbound}
    \lim_{\delta\downarrow 0}\frac{K_2(\delta,\alpha)}{K_1(\delta,\infty)}=0,
    \end{equation}
    where $K_1(\delta,\infty)$ is the constant in Proposition~\ref{mainparti}
    (recall \eqref{K1ident}).
    \end{enumerate}
    \end{prop}

\begin{proofsect}{Proof.}
(i) As in the proof of Proposition~\ref{mainparti}, we may insert the indicator on $\{B_T\geq 2\delta\}$
in the expectation on the left-hand side of \eqref{refine2} and add a factor of $1+o(1)$.

Introduce the following measurable subsets of $C^+_0$, respectively, $\Ccal^+$:
    \begin{eqnarray}
    G^\geq_{\delta,\alpha}&=&\{g\in C^+_0\colon g(\delta)=0,
    \max g\geq\alpha\delta^{-\frac 12}\},\label{GCdef}\\
    F^\geq_{\delta,\alpha}&=&\Bigl\{(y,f)\in\Ccal^+\colon y\geq 2\delta,
    \max_{[0,\delta]}f\geq\alpha\delta^{-\frac 12}\Bigr\}\label{Fdeltadef}.
    \end{eqnarray}
Note from \eqref{Ehatdef} and \eqref{Fhatdefgq} that
\begin{equation}
\begin{array}{ll}
&\eve(\delta;T)\cap\evf^{\,\geq} (\delta,\alpha;T)\cap\{B_T\geq 2\delta\}
= \{L(T,-\,\cdot)\in G_{\delta}\} \cap \\[0.2cm]
&\Big(
\{L(T,B_T+\cdot)\in G^\geq_{\delta,\alpha}\}\cup
\Big[
\{(B_T,L(T,B_T-\cdot)|_{[0,B_T]}) \in F^\geq_{\delta,\alpha}\} \cap
\{L(T,B_T+\cdot) \in G_\delta\}
\Big]
\Big)
\end{array}
\end{equation}
with $G_\delta=\{g\in C^+_0\colon g(\delta)=0\}$.

Pick $a\in\R$ such that $\mu+\rho(a)=0$, i.e., $a=\rho^{-1}(-\mu)<a^{**}$.
Apply Proposition~\ref{repr} twice for $G^-=G_\delta$ and the two choices:
(1) $F=F^\geq_{\delta,\alpha}$, $G^+=G_{\delta}$; (2) $F=\Ccal^+$,
$G^+=G^\geq_{\delta,\alpha}$. Sum the two resulting equations, to obtain
    \begin{equation}
    \label{refine2a}
    \begin{aligned}
    \mbox{}&\mbox{l.h.s.~of \eqref{refine2}}=(1+o(1))\int_0^\infty\d t_1\int_0^\infty\d t_2
    \1_{\{t_1+t_2\leq T\}}e^{a(t_1+t_2)}\\
    &\times\widehat\E^a\Bigl(\Bigl[
    \frac{w_{G_{\delta}}(X_0,t_1)}{x_a(X_0)}
    \1_{\{A^{-1}(T-t_1-t_2)\geq 2\delta\}}
    \1_{\bigl\{\max_{[0,\delta]}X\geq\alpha\delta^{-\frac 12}\bigr\}}\\[0.2cm]
    &\qquad\qquad +\frac{w_{G^\geq_{\delta,\alpha}}(X_0,t_1)}{x_a(X_0)}\Bigr]
    \frac{w_{G_{\delta}}(Y_{T-t_1-t_2},t_2)}{x_a(Y_{T-t_1-t_2})}\Bigr).
    \end{aligned}
    \end{equation}
In the same way as in the proof of Proposition~\ref{mainparti}, we obtain that
(recall (\ref{K1ident}--\ref{ydeltadef}))
    \begin{equation}
    \label{refine2b}
    \lim_{T\to\infty}\bigl(\mbox{r.h.s.~of \eqref{refine2a}}\bigr)
    = K_2(\delta,\alpha)
    \end{equation}
with
    \begin{equation}
    \label{K2def}
    K_2(\delta,\alpha) =
    \Bigl[\widehat\E^a_h\Bigl(\frac{y_a^{\smallsup{\delta}}(X_0)}{x_a(X_0)}
    \1_{\bigl\{\max_{[0,\delta]}X\geq\alpha\delta^{-\frac 12}\bigr\}}\Bigr)
    + \langle x_a,y_a^{\smallsup{\delta,\alpha}}\rangle\Bigr]
    \frac 1{\rho'(a)}\langle x_a, y_a^{\smallsup{\delta}}\rangle_\circ,
    \end{equation}
where $y_a^{\smallsup{\delta}}$ is defined in \eqref{ydeltadef} and
$y_{a}^{\smallsup{\delta,\alpha}}$ is defined as (recall \eqref{Wdef})
    \begin{equation}
    \label{ydeltaalpha}
    y_{a}^{\smallsup{\delta,\alpha}}(h)=\int_0^\infty\d t\,e^{at}
    w_{G^\geq_{\delta,\alpha}}(h,t)
    =\E_h^\star\Bigl(e^{\int_0^\infty[aX_v^\star- (X_v^\star)^2]\,\d v}
    \1_{\{X^\star_\delta=0\}}
    \1_{\bigl\{\max_{[0,\delta]}X^\star\geq\alpha\delta^{-\frac 12}\bigr\}}\Bigr).
    \end{equation}
The right-hand side of \eqref{K2def} is a strictly positive finite number.

\medskip\noindent
(ii) Fix $\alpha\in(0,\infty)$. From \eqref{K1ident} and \eqref{refine2b} we see that
$K_2(\delta,\alpha)/K_1(\delta,\infty)=K^{\smallsup{1}}(\delta,\alpha)
+K^{\smallsup{2}}(\delta,\alpha)$ with
    \begin{equation}
    \label{Ksdef}
    K^{\smallsup{1}}(\delta,\alpha)
    =\frac{\int_0^\infty\d h\, x_a(h)y_a^{\smallsup{\delta}}(h)
    \widehat\P^a_h\bigl(\max_{[0,\delta]}X\geq\alpha\delta^{-\frac 12}\bigr)}
    {\langle x_a,y_a^{\smallsup{\delta}}\rangle},
    \qquad K^{\smallsup{2}}(\delta,\alpha)
    =\frac{\langle x_a,y_a^{\smallsup{\delta,\alpha}}\rangle}
    {\langle x_a,y_a^{\smallsup{\delta}}\rangle}.
    \end{equation}
To prove \eqref{Kbound}, we need the following lemma.

    \medskip
    \begin{lemma}
    \label{prop-maxbd}
    Fix $a<a^{**}$ and $\alpha\in(0,\infty)$. Then:
    \begin{enumerate}
    \item[(i)]
    There exists $d=d(\alpha)>0$ such that, for any $R>0$ and any $\delta>0$
    sufficiently small,
    \begin{eqnarray}
    \sup_{h\in[0,R]}\widehat {\mathbb P}_h^a\Big(\max_{[0,\delta]}X
    \geq\alpha\delta^{-\frac 12}\Big)
    &\leq& c e^{-d \delta^{-\frac 14}}e^{c\sqrt R},\label{refinebound1}\\
    \sup_{h\in[0,R]}\frac{y_a^{\smallsup{\delta,\alpha}}(h)}{y_a(h)}&\leq&
    c e^{-d \delta^{-\frac 14}}e^{c\sqrt R}.\label{refinebound2}
    \end{eqnarray}
    \item[(ii)]
    For any $\delta>0$ sufficiently small,
    \begin{equation}
    \label{refinebound3}
    \inf_{h\in[0,\delta]}y_a^{\smallsup{\delta}}(h)\geq c.
    \end{equation}
    \end{enumerate}
    \end{lemma}

\begin{proofsect}{Proof.}
The proof is deferred to Section~\ref{sec-proof2}.
\end{proofsect}
\qed

\noindent
We use Lemma~\ref{prop-maxbd} to show that
\begin{equation}
\label{Klims}
\lim_{\delta\downarrow 0} K^{\smallsup{1}}(\delta,\alpha)
= \lim_{\delta\downarrow 0} K^{\smallsup{2}}(\delta,\alpha) = 0,
\end{equation}
which yields \eqref{Kbound}.

First note that, with the help of \eqref{refinebound3}, the common denominator
in \eqref{Ksdef} may be estimated from below by
    \begin{equation}
    \label{lowerdenom}
    \langle x_a,y_a^{\smallsup{\delta}}\rangle
    \geq \int_0^\delta \d h\, x_a(h) y_a^{\smallsup{\delta}}(h)
    \geq c \int_0^\delta \d h\, x_a(h) \geq c\delta,
    \end{equation}
where we use that $x_a$ is bounded away from zero on $[0,\delta]$.

In order to estimate the numerator of $K^{\smallsup{1}}(\delta,\alpha)$
from above, we split the integral in the numerator into two parts: $h\leq R$
and $h>R$. In the integral over $h\leq R$, estimate $y_a^{\smallsup{\delta}}\leq y_a$
and use \eqref{refinebound1}, to get the upper bound $ce^{-d\delta^{-\frac 14}}
e^{c\sqrt{R}}$. In the integral over $h>R$, estimate $y_a^{\smallsup{\delta}}\leq y_a$,
estimate the probability against one and use \eqref{xasyrough} and \eqref{yasyrough},
to get the upper bound $ce^{-cR^{\frac 32}}$. Pick $R$ such that $c\sqrt{R}=\frac d2
\delta^{-\frac 14}$, to obtain that the numerator of $K^{\smallsup{1}}(\delta,\alpha)$
is at most $ce^{-\frac d2 \delta^{-\frac 14}}$.

In the same way we show, with the help of \eqref{refinebound2}, that the numerator
of $K^{\smallsup{2}}(\delta,\alpha)$ in \eqref{Ksdef} is at most $ce^{-\frac d2
\delta^{-\frac 14}}$. Now combine the two estimates with \eqref{lowerdenom} to
obtain \eqref{Klims}.
\end{proofsect}
\qed

\section{Addendum 2: Proof of Lemma~\ref{prop-maxbd}}
\label{sec-proof2}

\setcounter{step}{0}
\noindent
We will need the
following asymptotics for $x_a$ and $y_a$, which are refinements of \eqref{xasyrough}
and \eqref{yasyrough}, respectively.

    \medskip
    \bes
    \label{xyasy}
    For $a<a^{**}$,
    \begin{equation}
    \lim_{h\to\infty} \frac{1}{\sqrt{h}}\log\bigl[e^{\frac{\sqrt{2}}{3}
    h^{\frac32}} x_{a}(h)\bigr]
    =\lim_{h\to\infty} \frac{1}{\sqrt{h}}\log\bigl[e^{\frac{\sqrt{2}}{3}
    h^{\frac32}} y_{a}(h)\bigr]
    =\frac{a}{\sqrt{2}}.
    \end{equation}
    \es

\begin{proofsect}{Proof.}
The statement for $y_a$ is well-known, and follows from \eqref{yident}
together with the asymptotics
of the Airy function given by (see \cite[p.~43]{Erde56})
    \begin{equation}
    \label{Airyasy}
    {\rm Ai}(h)=\frac{1}{2\pi h^{\frac 14}} e^{-\frac{2}{3}h^{\frac 32}}[1+o(1)],
    \qquad h\rightarrow \infty.
    \end{equation}
To prove the statement for $x_a$, use \cite[Theorem 2.1, pp.~143--144]{CL55}.
To this end, define
    \begin{equation}
    \label{zetadef}
    \zeta_1(h)= x_a(h^2), \qquad \zeta_2(h)= h^{-2}\zeta_1'(h).
    \end{equation}
Then the eigenvalue equation $\Kcal^ax_a=\rho(a)x_a$ (recall \eqref{opdef})
can be written as (see also \cite[equation (5.3)]{CL55})
    \begin{equation}
    \label{w'eq*}
    \zeta'(h) = h^2 B(h) \zeta(h),
    \end{equation}
with
    \begin{equation}
    \zeta(h)=\left(\begin{array}{ll}
                \zeta_1(h)\\
                \zeta_2(h)
                \end{array}
             \right),\qquad
    B(h)=\left(\begin{array}{cc}
                    0                                           &1\\
                    2 - \frac{2a}{h^2} + \frac{2\rho(a)}{h^4}   &-\frac{3}{h^3}
               \end{array}
         \right).
    \end{equation}
Note that $B(h)=\sum_{n=0}^\infty h^{-n}B^{\smallsup{n}}$ ($B^{\smallsup{0}}\neq 0$)
is a convergent power series in $h^{-1}$, with $B^{\smallsup{0}}$ having eigenvalues
$\lambda_{1,2}=\pm \sqrt{2}$. Therefore \eqref{w'eq*} has formal solutions of the form
    \begin{equation}
    Z(h)=P(h) h^R e^{Q(h)},
    \end{equation}
where the columns of the matrix $Z$ are the two linearly independent solutions to
\eqref{w'eq*}, $P(h)= \sum_{n=0}^{\infty} h^{-n} P^{\smallsup{n}}$ ($\det(P^{\smallsup{0}})
\neq 0$) is a formal power series in $h^{-1}$, $R$ is a complex diagonal matrix,
and $Q(h)= \frac{1}{3}h^3 Q^{\smallsup{0}} + \frac{1}{2}h^2Q^{\smallsup{1}}+
h Q^{\smallsup{2}}$ is a matrix polynomial with $Q^{\smallsup{0}}$,
$Q^{\smallsup{1}}$, $Q^{\smallsup{2}}$ diagonal. In our case,
    \begin{equation}
    Q^{\smallsup{0}}={\rm diag}\{-\sqrt{2},+\sqrt{2}\},\qquad
    Q^{\smallsup{1}}=0,\qquad
    Q^{\smallsup{2}}={\rm diag}\bigl\{\textstyle{\frac{a}{\sqrt{2}},-\frac{a}{\sqrt{2}}}\bigr\}.
    \end{equation}
{}From the proof of \cite[Theorem 2.1]{CL55} it follows that $P(h),R,Q(h)$ can be
chosen to be real, because $B(h),\lambda_{1},\lambda_2$ are real. On
\cite[p.~151]{CL55} there is the further remark that for every formal solution
there exists an actual solution with the same asymptotics.

We need the solution that is in $L^2[0,\infty)$. By construction, we compute,
for $R={\rm diag}\{r_1,r_2\}$ (with $r_1,r_2$ some functions of $a$),
    \begin{equation}
    h^R e^{\frac13 h^3Q^{\smallsup{0}}
    + \frac12 h^{2}Q^{\smallsup{1}}
    + h Q^{\smallsup{2}}}=\left(\begin{array}{cc}
    h^{r_1}e^{-\frac{\sqrt{2}}{3} h^{3}+\frac{a}{\sqrt{2}}h}      &0\\
    0  &h^{r_2}e^{\frac{\sqrt{2}}{3} h^{3}-\frac{a}{\sqrt{2}}h}
    \end{array}
    \right).
    \end{equation}
Therefore the solution in $L^2[0,\infty)$ must be
    \begin{equation}
    \zeta(h)=h^{r_1}e^{-\frac{\sqrt{2}}{3} h^{3}+\frac{a}{\sqrt{2}}h}
    \sum_{n=0}^{\infty} h^{-n} \left(\begin{array}{ll}
    P^{\smallsup{n}}_{11}\\
    P^{\smallsup{n}}_{21}\end{array}
    \right),
    \end{equation}
where $P^{\smallsup{n}}_{ij}$ denotes the element in the $i$-th row and the
$j$-th column of the matrix $P^{\smallsup{n}}$. Now return to \eqref{zetadef}
to read off the claim.
\end{proofsect}
\qed

Pick $a'$ such that $a<a'<a^{**}$ and define (recall \eqref{Ddef})
    \begin{equation}
    \label{Mdef}
    M_t=\frac{D_t^{\smallsup{a'}}}{D_t^{\smallsup{a}}}
    =\frac{x_a(X_0)}{x_{a'}(X_0)}\frac{x_{a'}(X_t)}{x_a(X_t)}
    e^{-t[\rho(a')-\rho(a)]}e^{(a'-a)\int_0^tX_v\,\d v},\qquad t\geq 0.
    \end{equation}

    \medskip
    \begin{step}
    \label{martingale} For any $h\ge 0$,
    $(M_t)_{t\geq 0}$ is a martingale under $\widehat\P_h^a$.
    \end{step}

\begin{proofsect}{Proof.}
Fix $0\leq s<t$. If $\phi_s$ denotes the time-shift by $s\geq0$ (i.e.,
$(\phi_s\circ f)\bigl((X_t)_{t\geq0}\bigr)=f\bigl((X_{s+t})_{t\geq 0}\bigr)$
for any bounded and measurable function $f$), then it is clear that $M_t=M_s
(\phi_s\circ M_{t-s})$. Hence, using the Markov property at time $s$, we see that,
for any $h\geq 0$,
    \begin{equation}
    \widehat \E_h^a(M_t\,|\,M_s)
    =M_s\widehat \E_h^a(\phi_s\circ M_{t-s}\,|\,M_s)
    =M_s \widehat\E_h^a\bigl(\widehat\E_{X_s}(M_{t-s})\bigr).
    \end{equation}
Now use that, for any $x\geq 0$, according to the construction of the transformed
process in (\ref{Ddef}--\ref{transdens}),
    \begin{equation}
    \label{mart1}
    \widehat \E_x^a(M_{t-s})=\E_x\bigl(D_{t-s}^{\smallsup{a}}M_{t-s}\bigr)
    =\E_x\bigl(D_{t-s}^{\smallsup{a'}}\bigr)=1.
    \end{equation}
\end{proofsect}
\qed

    \medskip
    \begin{step}
    \label{concl}
    Proof of \eqref{refinebound1}.
    \es

\begin{proofsect}{Proof.}
Use Step~\ref{martingale}, Doob's martingale inequality and \eqref{mart1}, to obtain
    \begin{equation}
    \label{Doob}
    \widehat\P_h^a\Bigl(\max_{[0,\delta]}M\geq K\Bigr)\leq \frac 1K \max_{t\in [0,\delta]}
    \widehat\E_h^a(M_t)=\frac 1K,\qquad h\geq 0,\, K>0.
    \end{equation}
Next note that by Step \ref{xyasy}, for any $R>0$,
    \begin{equation}
    \inf_{[0,R]}\frac{x_a}{x_{a'}} \geq e^{-c\sqrt{R}}.
    \end{equation}
Substitute this into \eqref{Mdef}, to get
    \begin{equation}
    \label{Mtest}
    M_t\geq c\frac{x_{a'}(X_t)}{x_a(X_t)}e^{-c\sqrt{R}}\quad\widehat\P_h^a\mbox{-a.s.},
    \,0\leq t\leq 1,\,0\leq h\leq R.
    \end{equation}
Pick $g_a\colon[0,\infty)\to(0,\infty)$ to be the largest increasing function not
exceeding $x_{a'}/x_a$ anywhere on $[0,\infty)$. Then, by \eqref{Mtest},
$M_t\geq c g_a(X_t) e^{-c\sqrt{R}}$ $\widehat\P_h^a$-a.s., $0\leq t\leq 1$,
$0\leq h\leq R$. Now use \eqref{Doob} to estimate, for $0\leq h\leq R$,
    \begin{eqnarray}
    \widehat \P_h^a\Bigl(\max_{[0,\delta]}X\geq \alpha\delta^{-\frac 12}\Bigr)
    &=&\widehat \P_h^a\Bigl(\max_{t\in[0,\delta]}c g_a(X_t)
    \geq c g_a\bigl(\alpha\delta^{-\frac 12}\bigr)\Bigr)\nonumber\\
    &\leq& \widehat \P_h^a\Bigl(\max_{t\in[0,\delta]}M_t\geq
    c g_a\bigl(\alpha\delta^{-\frac 12}\bigr)e^{-c \sqrt{R}}\Bigr)\nonumber\\
    &\leq&\frac{1}{cg_a\bigl(\alpha\delta^{-\frac 12}\bigr)}e^{c\sqrt{R}}.
    \end{eqnarray}
By Step~\ref{xyasy}, it is possible to pick $g_a$ such that
    \begin{equation}
    g_a(h) \geq e^{-c \sqrt h},\qquad h\to\infty.
    \end{equation}
This implies the bound in \eqref{refinebound1} with $d(\alpha)=\sqrt{\alpha}$.
\end{proofsect}
\qed

    \medskip
    \begin{step}
    Proof of \eqref{refinebound2}.
    \es

\begin{proofsect}{Proof.}
Fix $a<a^{**}$. Define
    \begin{equation}
    D^{\smallsup{a,\star}}_t=\frac {y_a(X_t^{\star})}{y_a(X_0^{\star})}
    e^{\int_0^t[aX_v^{\star}-(X_v^{\star})^2]\,\d v}, \qquad t\geq 0.
    \end{equation}
Then it is easy to check (see \cite[Section~VIII.3]{RY}) that
$(D^{\smallsup{a,\star}}_t)_{t\geq 0}$ is a martingale under $\P_h^\star$ for
any $h\geq 0$ ($y_a$ is a strictly positive solution to the differential equation
$2y_a''(h)=(h-a)y_a(h)$ on $[0,\infty)$; recall \eqref{Airydiff} and \eqref{yident}).
Hence, analogously to \eqref{transdens}, we may construct a transformed process
via a Girsanov transformation by taking $D^{\smallsup{a,\star}}_t$ formally as
a density with respect to BESQ$^0$. Denote by $\widehat\P^{a,\star}_h$
and $\widehat\E^{a,\star}_h$ probability and expectation with respect to this
transformed process starting at $h\ge 0$.

Recall that $y_a(0)=1$. We have the following representation for the function
$y_a^{\smallsup{\delta,\alpha}}$ (recall \eqref{ydeltaalpha}):
   \begin{equation}
   \label{yalphaident}
   y_a^{\smallsup{\delta,\alpha}}(h)
   = y_a(h)\widehat\P_h^{a,\star}\bigl(X^\star_\delta=0,
   \max X^\star\geq\alpha\delta^{-\frac 12}\bigr).
   \end{equation}
The proof of \eqref{refinebound2} is now analogous to Steps~\ref{martingale}--\ref{concl}.
Indeed, use \eqref{yalphaident}, drop the restriction $X^\star_\delta=0$, and proceed
analogously. Step~\ref{xyasy} provides the necessary asymptotic bounds for $y_a$ and
$y_{a'}$, provided that $a<a'<a^{**}$.
\end{proofsect}
\qed

    \medskip
    \begin{step}
    Proof of \eqref{refinebound3}.
    \es

\begin{proofsect}{Proof.}
We return to the right-hand side of \eqref{ydeltadef} and obtain a lower bound
by inserting the indicator of the event $\{\max X^\star\leq 2\delta\}$. On this
event, we may estimate the exponential from below by $c$. Hence, for $0\leq h\leq \delta$,
    \begin{equation}
    y_a^{\smallsup{\delta}}(h) \geq
    c\P^\star_h\bigl(\max X^\star\leq 2\delta,X^\star_\delta=0\bigr)
    = c\left[\P^\star_h(X_\delta^\star=0)
    -\P^\star_h\bigl(\max X^\star>2\delta,X^\star_\delta=0\bigr)\right].
    \end{equation}
Using the Markov property at the first time the BESQ$^0$ hits $2\delta$, we see
that the latter probability is at most $\P^\star_{2\delta}(X^\star_\delta=0)$. Since
the first probability is decreasing in $h$, we therefore have the bound
    \begin{equation}
    y_a^{\smallsup{\delta}}(h) \geq c\bigl(\P^\star_\delta(X^\star_\delta=0)
    -\P^\star_{2\delta}(X^\star_\delta=0)\bigr).
    \end{equation}
Now use that $\P^\star_{h}(X^\star_\delta=0)=e^{-h/2\delta}$ for any $h,\delta\geq 0$
(see \cite[Corollary XI(1.4)]{RY}), to complete the proof.
\end{proofsect}
\qed

%%% REFERENCES %%%

\end{document}